\documentclass{AIMS}
\usepackage{amsfonts}
\usepackage{amssymb,amsmath}
\usepackage{graphicx}
\usepackage{amsmath}
  \usepackage{paralist}
  \usepackage{graphics} 
  \usepackage{epsfig} 
 \usepackage[colorlinks=true]{hyperref}
\hypersetup{urlcolor=blue, citecolor=red}

\def\currentvolume{X}
 \def\currentissue{X}
  \def\currentyear{200X}
   \def\currentmonth{XX}
    \def\ppages{X--XX}

 \newcommand{\calP}{\mathcal{P}}
 \newcommand{\ext}{\rm ext}
 \newcommand{\sta}{\rm sta}
 \newcommand{\PP}{\Pi}
 \newcommand{\bx}{\boldsymbol{x}}
 \newcommand{\la}{\langle}
 \newcommand{\ra}{\rangle}
 \newcommand{\calX}{\mathcal{X}}
 \newcommand{\real}{\mathbb{R}}
 \newcommand{\Lam}{\Lambda}
 \newcommand{\calV}{\mathcal{V}}
 \newcommand{\VV}{V}
 \newcommand{\eb}{\begin{equation}}
  \newcommand{\ee}{\end{equation}}
  \newcommand{\bvsig}{\boldsymbol{\varsigma}}
  \newcommand{\bveps}{\mbox{\boldmath$\xi$}}
  \newcommand{\half}{\frac{1}{2}}
  \newcommand{\UU}{U}
  \newcommand{\bA}{\mathbf{A}}
  
  \newcommand{\calS}{\mathcal{S}}
  \newcommand{\calE}{\mathcal{E}}
  \newcommand{\calPd}{\calP^d}
  \newcommand{\barbx}{\bar{\bx}}
  \newcommand{\barbvsig}{\bar{\bvsig}}
  \newcommand{\barvsig}{\bar{\bvsig}}
  \newcommand{\PPd}{\PP^d}
  
  \newcommand{\veps}{\xi}
  \newcommand{\vsig}{\varsigma}
  \newcommand{\bff}{\boldsymbol{f}}
  \newcommand{\LL}{\boldsymbol{L}}
  \newcommand{\WW}{{W}}
  \newcommand{\bxi}{{\mbox{\boldmath$\xi$}}}
  
  \newcommand{\Diag}{{\mbox{Diag }}}
  \newcommand\barxi{\boldsymbol{\bar{\xi}}}
  \newcommand{\rank}{\rm rank}

\newtheorem{theorem}{Theorem}[section]

\newtheorem{lemma}[theorem]{Lemma}

\theoremstyle{definition}

\newtheorem{remark}{Remark}

  \newcommand{\bG}{{\bf G}}
  \newcommand{\G}{G\^{a}teaux}

\title[Triality Theory for Quartic Polynomial Optimization]
      {On the Triality Theory for a Quartic Polynomial Optimization Problem}

\author[D. Gao and C.Z. Wu]{}


 \email{d.gao@ballarat.edu.au}
 \email{changzhiwu@yahoo.com}

\begin{document}
\maketitle

\centerline{\scshape David Yang Gao }
\medskip
{\footnotesize
 \centerline{School of Science, Information Technology and Engineering,}
   \centerline{University of Ballarat, Victoria 3353, Australia.}
} 

\medskip

\centerline{\scshape Changzhi Wu}
\medskip
{\footnotesize
 \centerline{School of Science, Information Technology and Engineering,}
   \centerline{University of Ballarat, Victoria 3353, Australia.}
   \centerline{School of Mathematics, Chongqing Normal University,}
   \centerline{Shapingba, Chongqing, 400047, China}

\bigskip

 \centerline{(Communicated by Kok Lay Teo)}

\begin{abstract}
This paper presents a  detailed  proof of the triality theorem for
a class of fourth-order polynomial optimization problems. The
method is based on linear algebra but it solves an open problem on
the double-min duality left in 2003. Results show that the triality  theory
holds strongly  in a tri-duality form  if the
primal problem and its canonical dual have the same dimension;
otherwise, both the canonical min-max duality and the double-max
duality still hold strongly, but the double-min duality holds
weakly in a symmetrical form.
 Four numerical examples are  presented to illustrate
 that  this 
theory can be used to identify not only the  global minimum, but
also  the largest  local minimum  and local maximum.
\end{abstract}

\section{Introduction and Motivation}

The concepts of {\em triality} and {\em tri-duality} were
originally proposed in nonconvex mechanics \cite{gao-tri96,
gao-amr}. Mathematical theory of triality  in its standard format
is  composed  of three types of dualities:  a canonical min-max
duality 
 and a pair of     double-min  and  double-max dualities.
  The canonical min-max duality
  provides a sufficient condition for global minimum,
  while the double-min and double max dualities can be used to
 identify  respectively the largest local minimum and local maximum.
The tri-duality is a strong form of the triality principle
\cite{GaoBook}. Together with a {\em canonical dual
transformation}
 and a {\em complementary-dual principle,} they comprise a versatile
{\em canonical duality theory,}
 which can be used not only for solving a large class of challenging problems
 in nonconvex/nonsmooth analysis  and continuous/discrete optimization  \cite{GaoBook,gao-jogo00},
 but also for modeling complex systems and understanding multi-scale
 phenomena within a unified framework  \cite{gao-tri96,gao-amr,GaoBook}
(see also the review articles
\cite{gao-amma03,gao-cace09,Gao-Sherali-AMMA09}).

For example, in the recent work by Gao and Ogden
\cite{gao-ogden-qjmam08} on nonconvex variational/boundary value
problems, it was discovered that   both the global and local
minimizers are usually
 nonsmooth functions and cannot be determined easily  by traditional  Newton-type numerical methods.
However, by the canonical dual transformation, the nonlinear
differential equation is equivalent to an algebraic equation,
which can be solved analytically to obtain all solutions. Both
 global minimizer and local extrema were identified by the triality theory, which revealed   some  interesting
phenomena in phase transitions.

The triality theory has attracted much attention recently in
duality ways: Successful applications in multi-disciplinary fields
of  mathematics, engineering and sciences show  that this theory
is not only  useful and versatile, but also  beautiful in  its
mathematical  format and rich in connotation of physics,
 which reveals a unified intrinsic  duality pattern in complex systems;
 On the other hand,
  a large number of
    ``counterexamples"  have been presented 
     in  several papers since 2010. Unfortunately, most of these counterexamples are either
     fundamentally wrong (see \cite{svz,VZ-AA}), or
      repeatedly address an open problem left by Gao in 2003 on the double-min duality \cite{gao-opt03,gao-amma03}.

The main goal of this paper is to
  solve this open problem left in 2003.
The next section will present a brief review
  and the open problem in the triality theory.
   In Section 3, the triality theory is proved in  its strong form
  as it was   originally discovered.
  Section 4  shows that both the canonical min-max   and the
  double-max dualities hold strongly in general,
    but the double-min duality holds weakly in a symmetrical form.
  Applications are illustrated in Section 5, where a  linear perturbation
   method is used for solving certain critical problems.
  The paper ended by an Appendix and a section of concluding remarks.

\section{Canonical Duality Theory: A Brief Review and an Open Problem}
Let us   begin  with the  general global extremum  problem
\begin{equation}
(\calP): \;\; \ext \left\{ \PP\left( \bx\right) =W\left(
\bx\right) +\frac{1}{2} \la \bx , \mathbf{A}\bx \ra -  \la \bx ,
\boldsymbol{f}  \ra \; | \; \bx \in \calX_a \right\} , \label{pp}
\end{equation}
where  $\calX_a \subset \real^n$ is an open set,
 $\bx= \{ x_i \} \in
\mathbb{R} ^{n}$  is a decision vector, $\mathbf{A}=\left\{
A_{ij}\right\} \in \mathbb{R}^{n\times n}$ is a given symmetric
matrix, $\boldsymbol{f} =\{ f_i \} \in \mathbb{R}^{n}$ is a given
vector, and $\la * , * \ra $ denotes a bilinear form on $\real^n
\times \real^n$; the function $W: \calX_a \rightarrow \real $ is
assumed to be nonconvex and differentiable (it is allowed to be
nonsmooth and sub-differentiable for  constrained problems).
 The  notation $\ext \{ * \}$ stands for finding global extremal of the function given in $\{ * \}$.

  In this paper, we are interested only in three types of  global extrema:
  the global minimum and a pair of the largest local minimum and local maximum.
Therefore,  the nonconvex term  $W(\bx)$ in (\ref{pp}) is  assumed
to satisfy the objectivity condition \footnote{The concept of
objectivity in science means that qualitative and quantitative
descriptions of physical phenomena remain unchanged when the
phenomena are observed under a variety of conditions. That is, the
objective function should be independent with the choice of the
coordinate systems. In continuum mechanics, the objectivity is
also regarded as the {\em principle of frame-indifference. } See
Chapter 6 in \cite{GaoBook} for mathematical definitions of the
objectivity and geometric nonlinearity in differential geometry
and finite deformation field theory. Detailed discussion of
objectivity in global optimization will be given in another paper
\cite{gss-11}.}, i.e., there exists a  {\em (geometrically)
nonlinear  mapping} $\Lam:  \calX_a \rightarrow \calV \subset
\real^m$ and a canonical function $\VV: \calV \subset \real^m
\rightarrow \real$ such that
\begin{equation*}
W\left( \bx\right) = V\left( \Lambda \left( \bx\right) \right)
\;\; \forall \bx \in \calX_a.
\end{equation*}
According to \cite{GaoBook}, a real valued function $\VV:\calV
\rightarrow \real$ is said to be a canonical function on its
effective domain $\calV_a \subset \calV$ if its Legendre conjugate
$\VV^* :\calV^* \rightarrow \real$
 \eb
  \VV^*(\bvsig) = \sta \{ \langle \bveps ;
\bvsig \rangle  - \VV(\bveps) | \; \bveps \in \calV_{a} \} \ee
 is
uniquely defined on its effective domain $\calV^*_a \subset
\calV^*$ such that the canonical duality relations \eb \bvsig =
\nabla \VV(\bveps) \;\;  \Leftrightarrow \;\; \bveps = \nabla
\VV^*(\bvsig)  \;\; \Leftrightarrow \;\; \VV(\bveps ) +
\VV^*(\bvsig)  = \la \bveps ; \bvsig \ra
 \ee
hold on $\calV_a \times \calV^*_a$, where $\la * ;  * \ra $
represents a bilinear form which puts $\calV$ and $\calV^*$ in
duality.  The notation $\sta\{ *  \}$ stands for solving the
stationary point problem in $\{  *  \}$. By this one-to-one
canonical duality, the nonconvex function $ W(\bx) =
\VV(\Lam(\bx)) $  can be replaced by
  $\la \Lam(\bx) ; \bvsig \ra - \VV^*(\bvsig)$
such that the nonconvex function $\PP(\bx)$ in (\ref{pp}) can be
written  as \eb \Xi(\bx, \bvsig) = \la \Lam(\bx) ; \bvsig \ra -
\VV^*(\bvsig) + \half \la \bx , \mathbf{A} \bx \ra - \la \bx
,\boldsymbol{f} \ra , \ee which is the so-called {\em total
complementary (energy)  function} introduced by Gao and Strang in
1989. By using this total complementary function,
  the canonical dual function $\PP^d:\calV^*_a \rightarrow \real$
can be formulated as
  \eb
   \PP^d(\bvsig) = \sta \{ \Xi(\bx, \bvsig)
| \;\; \forall \bx \in \calX_a\} = \UU^\Lam(\bvsig) -
\VV^*(\bvsig), \label{Canonical Dual}
 \ee
  where $\UU^\Lam:\calV^*_a \rightarrow \real$ is called the
$\Lam$-conjugate of $\UU(\bx) =  \la \bx , \boldsymbol{f} \ra -
\half \la  \bx , \bA \bx \ra $, defined by \cite{GaoBook} as
 \eb
  \UU^\Lam(\bvsig) = \sta \{
  \la \Lam(\bx) ; \bvsig \ra - \UU(\bx) | \; \bx \in \calX_a \}.
 \ee
Let $\calS_a \subset \calV^*_a$ be the feasible domain of
$\UU^\Lam(\bvsig)$; then the canonical dual problem is to solve
the   stationary point problem
 \eb
 (\calPd): \;\; \ext \{ \;
\PP^d(\bvsig) | \; \bvsig \in \calS_a \}. \ee

\begin{theorem}[Complementary-duality principle \cite{GaoBook}]
\label{CDP}
 Problem $(\calPd)$ is a canonical dual to $(\calP)$
  in the sense that if $(\barbx, \barvsig)$ is a critical
point of $\Xi(\bx, \bvsig)$, then $\barbx$  is a critical point of
$(\calP)$, $\barvsig$ is a critical point of $(\calPd)$, and \eb
\PP(\barbx) = \Xi(\barbx, \barvsig)  = \PPd( \barvsig). \ee
\end{theorem}

Theorem \ref{CDP} implies a perfect duality relation (i.e. no
duality gap) between the primal problem and its canonical
dual\footnote{The complementary-dual in physics  means
 perfect dual in optimization, i.e.,
  the canonical dual in Gao's work,
  which means no duality gap. Otherwise, any duality gap will violet the energy conservation law.
  Therefore,  each complementary-dual variational statement  in continuum mechanics
is usually refereed as a principle.}.
 The formulation
of   $\PPd(\bvsig)$ depends on the geometrical operator
$\Lam(\bx)$.
 In many applications,
the geometrical operator  $\Lam$ is usually a quadratic mapping
over a given field \cite{GaoBook}. In finite dimensional space,
this quadratic operator can be written  as  a vector-valued
function (see \cite{gao-jogo00}, page 150)
 \eb
  \Lam(\bx) = \left\{ \half \bx^T \mathbf{B}^k \bx \right\}_{k=1}^m : \;\;\calX_a \subset
\real^n \rightarrow \calV_a  \subset \real^m , \label{eq-quadLam}
\ee where $\mathbf{B}^k = \{ \boldsymbol{B}^k_{ij} \} \in
\real^{n\times n} $ is a symmetrical matrix for each $k = 1, 2,
\cdots, m$, and $\calV_a \subset \real^m$ is defined by
\[
\mathcal{V}_{a}=  \left\{ \bveps \in \real^m | \;\; \veps_k =
\half \bx^T \mathbf{B}^k \bx \;\; \forall \bx \in \calX_a, \;\; k
= 1, \dots, m \right\} .
\]
 In this case,
  the total complementary function has the form
\eb \Xi \left( \bx,\bvsig\right) = \frac{1}{2} \la \bx ,
\mathbf{G} \left( \bvsig\right) \bx  \ra - \VV^*(\bvsig) - \la
\bx, \boldsymbol{f} \ra ,   \label{Lag} \ee where $\mathbf{G}:
\real^m \rightarrow \real^{n\times n}$ is a matrix-valued function
defined by
 \eb
\mathbf{G}(\bvsig) = \mathbf{A} + \sum_{k=1}^m \vsig_k
\mathbf{B}^k   . \label{Gkesi} \ee The critical condition
$\nabla_{\bx}\Xi \left( \bx, \bvsig\right) =0$ leads to the
canonical equilibrium equation
\begin{equation}
\mathbf{G} \left( \bvsig\right) \bx=\boldsymbol{f} . \label{Equi}
\end{equation}
Clearly, for any given $\bvsig\in \calV^*_a,$ if the vector $
\boldsymbol{f}\in \mathcal{C}_{ol}\left( \mathbf{G} \left(
\bvsig\right) \right) $, where $\mathcal{C}_{ol}\left( \mathbf{G}
\right)$ stands for a  space spanned by the columns of
$\mathbf{G}$, the canonical equilibrium equation (\ref{Equi}) can
be solved analytically as\footnote{In this paper $\mathbf{G}^{-1}$
should be understood as the generalized inverse if $\det
\mathbf{G} = 0$ \cite{gao-jogo00}.} $\bx =
[\mathbf{G}(\bvsig)]^{-1} \bff $. Therefore, the canonical dual
feasible space $\mathcal{S}_{a}\subset \calV^*_a $ can be defined
as
\begin{equation*}
\mathcal{S}_{a}=\left\{ \bvsig\in \mathcal{V}_{a}^{\ast } \;|\;
\boldsymbol{f}\in \mathcal{C}_{ol}\left( \mathbf{G} \left(
\bvsig\right) \right)  \right\} ,
\end{equation*}
and on $\calS_a$  the canonical dual $\PP^d(\bvsig)$ is
well-defined  as \eb
 \PP^d(\bvsig) = - \half \la
\mathbf{G}(\bvsig)^{-1} \boldsymbol{f} ,  \boldsymbol{f} \ra -
\VV^*(\bvsig) . \label{eq-pid} \ee

\begin{theorem}[Analytic solution \cite{GaoBook}] \label{thm-solu}
   If
$\barvsig \in \calS_a$ is a critical  solution of $(\calPd)$,
 then
 \eb
 \barbx = \mathbf{G}(\barvsig)^{-1} \bff \label{eq-solu}
 \ee
is a critical  solution of $(\calP)$  and $\PP(\barbx) =
\PPd(\barvsig)$.

Conversely, if $\barbx$  is a critical solution of $(\calP)$, it
must be in the form of  (\ref{eq-solu})  for a certain critical
solution  $\barvsig$  of $(\calPd)$.
\end{theorem}

The canonical dual function $\Pi^d(\bvsig)$ for a general
quadratic operator $\Lam$  was first formulated in nonconvex
analysis, where  Theorem \ref{thm-solu} is called the pure
complementary energy principle, \cite{gao-mecc99}. In finite
deformation theory, this theorem
 solved an open problem left by Hellinger (1914) and Reissner
(1954)  (see \cite{li-cupta}). The analytical solution theorem has
been successfully applied for solving a class of nonconvex
problems in mathematical physics, including Einstein's special
relativity theory \cite{GaoBook},
  nonconvex mechanics and phase transitions in solids \cite{gao-ogden-qjmam08}.
 In global optimization, the primal solutions to
nonconvex minimization and integer programming problems are
usually located on the boundary of the feasible space.
 By Theorem \ref{thm-solu},
these solutions can be analytically determined by critical points
of the canonical dual function $\PPd(\bvsig)$ (see
\cite{fang-gaoetal07, gao-jimo07, gao-ruan-jogo10, g-r-s-JOGO10}).

In order to identify both global and local extrema  of the primal
and dual problems, we assume, without losing much generality, that
the canonical function $\VV:\calE_a \rightarrow \real$ is convex
and let
\begin{eqnarray}
 \mathcal{S}_{a}^{+} &=& \left\{ \bvsig\in
\mathcal{S}_{a}\;|\;\mathbf{G}\left( \bvsig\right)  \succeq
0\right\}  \label{SaPos},  \\
 \mathcal{S}_{a}^{-} &=& \left\{ \bvsig\in
\mathcal{S}_{a}\;|\;\mathbf{G}\left( \bvsig\right) \prec 0\right\}
\label{Sa-} ,
\end{eqnarray}
where $\mathbf{G}\left( \bvsig\right)  \succeq 0$ means that
$\mathbf{G}\left( \bvsig\right)$ is a positive semi-definite
matrix and $\mathbf{G}\left( \bvsig\right) \prec 0$ means that
$\mathbf{G}\left( \bvsig\right)$ is  negative definite.

\begin{theorem} [Triality Theorem \cite{gao-jogo00}]
Let $(\barbx,\barvsig)$
 be a critical point of $\Xi \left(
\bx,\bvsig\right)$.

If $ \bG(\barvsig) \succeq 0 $,
 then $\barvsig$ is a global
maximizer of Problem $ ( \mathcal{P}^{d} )$, the vector  $\barbx$ is
a global minimizer of Problem $\left( \mathcal{P} \right)$, and
the following canonical   min-max duality statement holds:
\begin{equation}
\min_{\bx\in \mathcal{X}_{a}}\PP\left( \bx\right) =\Xi \left(
\barbx, \barvsig\right) =\max_{ \boldsymbol{\varsigma }\in
\mathcal{S}_{a}^{+}}\PP^{d}\left( \bvsig \right) .  \label{Global0}
\end{equation}

If $ \bG(\barvsig ) \prec 0 $, then there exists a
neighborhood $\mathcal{X}_{o}\times \mathcal{S}_{o}\subset
\calX_a \times \mathcal{S}_{a}^{-}$ of $\left(
\barbx,\barvsig\right) $ for which we have either the double-min duality statement
\begin{equation}
\min_{\bx\in \mathcal{X}_{o}}\PP\left( \bx\right) =\Xi \left(
\barbx,\barvsig\right) =\min_{ \bvsig\in
\mathcal{S}_{o}}\PP^{d}\left( \bvsig\right) , \label{dmin}
\end{equation}
or the double-max duality statement
\begin{equation}
\max_{\bx\in \mathcal{X}_{o}}\PP\left( \bx\right) =\Xi \left(
\barbx,\barvsig\right) =\max_{ \bvsig\in
\mathcal{S}_{o}}\PP^{d}\left( \bvsig\right). \label{dmax}
\end{equation}

\end{theorem}

The triality theory provides actual
global extremum criteria for three types of solutions to the
nonconvex problem $(\calP)$:
  a  global minimizer $\barbx(\barbvsig)$ if $\barbvsig \in \calS^+_a$
  and a pair of the largest-valued local extrema. In other words,
   $\barbx(\barbvsig)$  is the largest-valued local maximizer if $\barbvsig \in \calS^-_a$ is
a  local maximizer;   $\barbx(\barbvsig)$ is the largest-valued local minimizer
    if $\barbvsig \in \calS^-_a$ is a local  minimizer.
This pair of largest local extrema plays a critical  role in nonconvex
analysis of post-bifurcation and phase transitions.

 \begin{remark}[
 Relation between Lagrangian Duality and  Canonical Duality]\hfill
 {\em

 The main difference between   the
  Lagrangian-type dualities (including the equivalent Fenchel-Moreau-Rockfellar dualities)
  and the canonical duality  is
  the  operator $\Lam:\calX_a \rightarrow \calV_a$.
  In fact, if $\Lam$ is linear, the primal problem $(\calP)$ is called geometrically linear
  in \cite{GaoBook} and  the
  total complementary function $\Xi(\bx, \bvsig)$ is simply  the well-known Lagrangian
  and is denoted as
  \eb
  \LL(\bx, \bvsig) = \la \Lam \bx ; \bvsig \ra - \VV^*(\bvsig) - F(\bx).
  \ee

  In convex (static) systems, $F(\bx) = \la \bx , \bff \ra $ is linear and
  $\LL(\bx, \bvsig)$ is a saddle function. Therefore, the well-known saddle min-max
  duality links a convex minimization problem $(\calP)$ to a concave maximization dual problem
  with linear constraint:
  \eb
  \max \left \{ \Pi^*(\bvsig) =  - \VV^*(\bvsig) | \;\; \Lam^* \bvsig = \bff, \;\; \bvsig \in \calV_a^* \right\}, \label{eq-pis}
  \ee
  where $\Lam^*$ is the conjugate operator of  $\Lam$ defined via
  $\la \Lam \bx ; \bvsig \ra = \la \bx , \Lam^* \bvsig \ra$.
  Using the Lagrange multiplier $\bx \in \calX_a$ to relax the equality constraint,
  the Lagrangian $\LL(\bx, \bvsig)$ is obtained.
  By the fact that the (canonical) duality in convex static systems is unique, the saddle min-max
  duality is   refereed as the {\em mono-duality}  in complex systems (see Chapter 1 in \cite{GaoBook}).

Since the linear operator $\Lam$ can not change the convexity of
$\WW(\bx) = \VV(\Lam\bx)$, the Lagrangian duality theory can be
used mainly  for convex problems. It is known that if $\WW(\bx)$
is nonconvex, then the Lagrangian duality as well as the related
Fenchel-Moreau-Rockafellar duality will produce the so-called
duality gap. Comparing the canonical dual function $\Pi^d(\bvsig)$
in (\ref{eq-pid}) with the Lagrangian dual function
$\Pi^*(\bvsig)$  in  (\ref{eq-pis}), we know that the duality gap
is $\half \la \mathbf{G}(\bvsig)^{-1} \boldsymbol{f} ,
\boldsymbol{f} \ra $.

The canonical duality theory is based on the (geometrically)
nonlinear mapping $\Lam:\calX_a \rightarrow \calV_a$ and the
canonical transformation $\WW(\bx) = \VV(\Lam(\bx))$.
 The total complementary function
$\Xi(\bx, \bvsig)$  is also known  as the   nonlinear or extended
Lagrangian and is denoted by $L(\bx, \bvsig)$ due to the geometric
nonlinearity of $\Lam(\bx)$ (see \cite{GaoBook, gao-amma03}).
Relations between the canonical duality and the classical
Lagrangian duality are discussed in \cite{g-r-s-JOGO10}.}
\end{remark}

\begin{remark}[Geometrical Nonlinearity and Complementary Gap Function]
{\em
 The canonical min-max duality statement (\ref{Global0}) was first
proposed by Gao and Strang in nonconvex/nonsmooth analysis and mechanics
 in 1989 \cite{GaoStrang89}, where $\Pi(\bx) = \WW(\bx) -
 F(\bx) $ is the so-called total potential energy with  $\WW$
representing  the internal (or stored) energy and $F$ the
external energy.
The {\em geometrical nonlinearity} is a standard terminology in finite deformation theory,
which  implies that the geometrical equation (or the  configuration-strain relation)
$\bxi = \Lam(\bx)$ is nonlinear.
By definition in physics,  a function $ F(\bx)$ is called the
 external energy means that its (sub-)differential must be the external force (or input)
 $\bff$.  Therefore, in Gao and Strang's work,
  the external energy should be a linear function(al) $ F(\bx) =   \la \bx, \bff \ra $
on its effective domain. In this
case, the matrix $\bG(\bvsig) $ is a Hessian of the so-called {\em
complementary gap function} (i.e. the Gao-Strang gap function
\cite{GaoStrang89})
\eb
  G_{ap} (\bx, \bvsig) =    \la  - \Lam_c
(\bx) ; \bvsig  \ra , \label{eq-gap}
\ee
 where $ \Lam_c(\bx) = - \half \bx^T {\bf B}^k \bx$ is called the complementary operator of a \G
differential $\Lam_t(\bx) \bx = \bx^T {\bf B}^k \bx$ of $\Lam(\bx)$
\cite{GaoStrang89}. Actually, in Gao and Strang's original work,
the canonical min-max duality statement holds in a general (weak)
condition, i.e.,   $ G_{ap} (\bx, \barbvsig)  \ge 0,\; \forall
\bx \in \calX_a$ in field theory (corresponding to the strong
condition $\bG(\barbvsig) \succeq 0$). The related canonical
duality theory has been generalized to nonconvex variational
analysis of a large deformation (von Karman) plate (where $\bA =
\Delta^2$ \cite{Yau-Gao}), nonconvex (chaotic)  dynamical systems
(where  $\bA = \Delta - \partial^2/\partial t^2$
\cite{gao-amma03}),  and general nonconvex constrained
problems in global optimization.
Since $F(\bx)$ in these general applications is the  quadratic function $- \half \la \bx , \mathbf{A} \bx \ra + \la \bx, \bff \ra$,
the Gao-Strang  gap function (\ref{eq-gap}) should be replaced
by the generalized form
$ G_{ap} (\bx, \bvsig) =   \half  \la  \bx , \bG(\bvsig) \bx \ra $ (see the review article by Gao and
Sherali \cite{Gao-Sherali-AMMA09}).
 This gap
function recovers the existing  duality gap in traditional
 duality theories and
provides a sufficient global optimality condition for general
nonconvex problems in both infinite and finite dimensional systems
(see review articles \cite{gao-amma03,Gao-Sherali-AMMA09}).
By the fact that the geometrical mapping $\Lam$ in Gao and Strang's work is a tensor-like
operator, it  has been realized recently that the popular semi-definite programming method is actually   a
special application (where   $\WW(\bx)$ is a quadratic function)
of the canonical min-max duality theory proposed in
1989 (see \cite{gao-ruan-jogo10,gao-ruan-pardalos}).

 In a recent paper by Voisei and Zalinescu \cite{VZ-AA}, they unfortunately
 misunderstood  some basic
 terminologies in continuum physics, such as
geometric nonlinearity, internal and external energies, and
present ``counterexamples" to the Gao-Strang theory based on certain  ``artificially chosen" operators
 $\Lam(\bx)$ and quadratic functions $F(\bx)$.
Whereas  in the stated contexts, the geometrical operator $\Lam$ should be  a  canonical measure
(Cauchy-Reimann type
finite deformation operator, see Chapter 6 in \cite{GaoBook})
and   the external energy $F(\bx)$  is typically a linear functional on its effective domain;
otherwise, its (sub)-differential will not be the external force.
 Interested readers are
refereed by \cite{gss-11} for further discussion. }
\end{remark}

\begin{remark}[Double-Min Duality and Open Problem]
{\em The double-min and double-max duality statements were
discovered simultaneously in a post-buckling analysis of large
deformed beam model  \cite{gao-tri96,gao-amr}  in 1996,
 where the finite
 strain measure $\Lam$ is a quadratic differential operator from a
2-D displacement field to a 2-D canonical strain field. Therefore,
the triality theory was first proposed  in its strong form,
 i.e. the so-called {\em tri-duality theory} (see the next section).
 Later on  when Gao was writing his duality book \cite{GaoBook}, he
 realized that this pair  of double-min and double-max dualities holds naturally
 in convex Hamilton systems.
 Accordingly,  a {\em bi-duality theorem} was proposed and proved
 for geometrically linear  systems (where $\Lam$  is a linear operator;
 see 
 Chapter 2 in \cite{GaoBook}).
Following this, the triality  theory was naturally generalized
  to
geometrically nonlinear systems (nonlinear $\Lam$; see Chapter 3
in \cite{GaoBook}) with applications  to global optimization
problems \cite{gao-jogo00}. However, it was discovered in 2003
 that if $ n \neq  m$ in the quadratic mapping
(\ref{eq-quadLam}),   the double-min duality statement needs
``certain additional constraints". For the sake of mathematical
rigor,
 the double-min duality was not included in
the triality theory and these additional constraints were left as
an open problem (see Remark 1 in \cite{gao-opt03},  also  Theorem
3 and its Remark in a review article by Gao \cite{gao-amma03}). By
the fact  that the double-max duality is always true,
the double-min duality was still included in the triality theory
 in the ``either-or"   form in many applications (see \cite{gao-cace09,g-r-s-JOGO10}).
 However, ignoring the open problem related to the ``certain additional constraints" on the double-min duality statement has led to some misleading results.}
 \end{remark}


  The goal of this paper is to solve this open problem
   by providing a simple proof of the triality theory based on  linear
algebra.
   To help understanding
the intrinsic characteristics of the original problem and its
canonical dual, we assume that the nonconvex objective function
$W(\bx)$ is a sum of  fourth-order canonical polynomials
 \begin{equation}
 \WW(\bx)  =\frac{1}{2}\sum\limits_{k=1}^{m} \beta _{k}\left(
\frac{1}{2}\bx^{T}\mathbf{B}^{k}\bx-  d^{k} \right) ^{2} ,
\label{P1}
\end{equation}
where   $\mathbf{B}^{k}=\left\{ B_{ij}^{k}\right\} \in
\mathbb{R}^{n\times n}$, $k=1,\cdots ,m,$ are all symmetric
matrices,   $\beta _{k}>0 $  and $d^{k}\in \mathbb{R},$
$k=1,\cdots ,m $ are given constants.
 This polynomial is actually a discretized form of the so-called
 {\em double-well potential}, first proposed by  van der Waals
in thermodynamics in 1895 (see \cite{venderWaals}), which  is the
mathematical model for natural phenomena of bifurcation and phase
transitions in
 biology, chemistry,  cosmology, continuum mechanics,
 material science, and quantum field theory, etc.
 (see \cite{gao-amr,Gao-Yu2008,Jaffe,Kibble}).
By using  the quadratic geometrical   operator $\Lam(\bx)$ given
by (\ref{eq-quadLam}),
  the canonical function
\begin{equation}
V\left( \bxi\right) =\frac{1}{2}(\bxi- {\boldsymbol{d}} )^{T}
\boldsymbol{\beta} (\bxi- {\boldsymbol{d}} )\; \label{VKesi}
\end{equation}
  and   its Legendre conjugate
\begin{equation}
V^{\ast }\left( \bvsig\right) = \frac{1}{2}\bvsig^{T}
\boldsymbol{\beta}^{-1}\bvsig+ \bvsig^{T}\boldsymbol{d}
\label{Conj}
\end{equation}
are  quadratic functions, where $\boldsymbol{\beta} = \Diag
(\beta^k)$ represents the diagonal matrix defined by the non-zero
vector $\{ \beta^k \}$.

In  the following discussions, we assume that all the critical
points of problem ($\calP$) are  non-singular, i.e., if
$\nabla\PP(\barbx) = 0$, then
 \eb
 \det \nabla^2\PP(\barbx) \neq 0. \label{non-deg}
 \ee
We will first prove that  if $n=m$, the triality theorem holds in
its strong form; otherwise,  the theorem holds in its weak form.
Three numerical examples are used to illustrate the effectiveness
and efficiency of the canonical duality theory.

\section{Strong Triality Theory for Quartic Polynomial Optimization: Tri-Duality Theorem }

We first consider the case   $m=n$. For simplicity, we assume that
$\beta _{k}=1$ in the following discussion (otherwise,
$\mathbf{B}^{k}$ can be replaced by $\sqrt{\beta
_{k}}\mathbf{B}^{k}$ and $ d^{k}$ is replaced by
$d^{k}/\sqrt{\beta _{k}} $ ). In this case,
 the problem (\ref{pp}) is denoted as problem ($\mathcal{P}$). Its canonical dual is
\begin{eqnarray}
 ( \mathcal{P}^{d} ) : & \ext  &  \left\{ \PP ^{d}\left(
\bvsig\right) =-\frac{1}{2}\boldsymbol{f}^{T}\left[ \mathbf{G}
\left( \bvsig\right) \right]
^{-1}\boldsymbol{f}-\frac{1}{2}\bvsig^{T}\bvsig-\bvsig^{T}
\boldsymbol{d}\;|\;\bvsig\in \mathcal{S}_{a} \subset \real^n
\right\} . \label{Dual}
\end{eqnarray}

\begin{theorem} [Tri-Duality Theorem] \label{STT} ~~~

Suppose that $m=n$, that the assumption (\ref{non-deg}) is
satisfied, that $\barvsig$
 is a critical point of  Problem
$ ( \mathcal{P} ^{d} ) $ and that $\barbx =\left[ \mathbf{G}
\left( \barvsig\right) \right] ^{-1}\boldsymbol{f}$.

If $\barvsig\in \mathcal{S}_{a}^{+},$ then $\barvsig$ is a global
maximizer of Problem $ ( \mathcal{P} ^{d} )$ in
$\mathcal{S}_{a}^{+}$ if and only if $\barbx$ is a global
minimizer of Problem $\left( \mathcal{P} \right)$, i.e.,
 the following canonical min-max statement holds:
\begin{equation}
\PP (\barbx)=\min_{\bx\in \mathbb{R} ^{n}}\PP \left( \bx\right)
\Longleftrightarrow  \max_{ \bvsig\in \mathcal{S}_{a}^{+}}\PP
^{d}\left( \bvsig \right) =  \PP ^{d}(\barvsig) . \label{Global}
\end{equation}

On the other hand, if $\barvsig\in \mathcal{S}_{a}^{-}$, then,
there exists a neighborhood $\mathcal{X}_{o}\times
\mathcal{S}_{o}\subset \mathbb{R} ^{n}\times \mathcal{S}_{a}^{-}$
of $\left( \barbx,\barvsig\right) $, such that either one of the
following two statements holds.

(A) The double-min duality statement
\begin{equation}
\PP (\barbx)=\min_{\bx\in \mathcal{X}_{o}}\PP \left( \bx\right) \;
\Longleftrightarrow  \;  \min_{ \bvsig\in \mathcal{S}_{o}}\PP
^{d}\left( \bvsig\right) = \PP ^{d}\left( \barvsig\right) ,
\label{double-min}
\end{equation}
 or (B) the double-max duality statement
\begin{equation}
\PP (\barbx) = \max_{\bx\in \mathcal{X}_{o}}\PP \left( \bx\right)
\; \Longleftrightarrow  \; \max_{ \bvsig\in \mathcal{S}_{o}}\PP
^{d}\left( \bvsig\right)  = \PP ^{d}\left( \barvsig\right) .
\label{double-max}
\end{equation}
\end{theorem}

\noindent \textbf{Proof. } If $\barvsig$ is a critical point of
the canonical dual problem $(\mathcal{P}^{d})$, the criticality
condition
\begin{equation}
\nabla \PP ^{d}\left( \boldsymbol{\varsigma }\right)
=\frac{1}{2}\left[
\begin{tabular}{c}
$\boldsymbol{f}^{T}\left[ \mathbf{G} \left( \bvsig\right) \right]
^{-1}\mathbf{B}^{1}\left[ \mathbf{G} \left( \bvsig\right) \right]
^{-1}
\boldsymbol{f}$ \\
$\cdots $ \\
$\boldsymbol{f}^{T}\left[\mathbf{G} \left( \bvsig\right) \right]
^{-1}\mathbf{B}^{n}\left[\mathbf{G} \left( \bvsig\right) \right]
^{-1}\boldsymbol{f}$
\end{tabular}
\ \right] -\boldsymbol{\varsigma }- \boldsymbol{d} =  0 \in
\real^n \label{GraDual}
\end{equation}
leads to $\barvsig = \Lam(\barbx)$. By the fact that $\nabla
\Pi(\barbx) = \mathbf{G}(\bar{\bvsig})\barbx - \bff  = 0 \in
\real^n$,
   it follows that   $\barbx =\left[
\mathbf{G} \left( \barvsig\right) \right] ^{-1}\boldsymbol{f}$ is
a critical point of Problem $(\mathcal{P} )$.

To prove the validity of the canonical min-max statement
(\ref{Global}), let $\barvsig$ be a critical point and
$\barvsig\in \mathcal{S}_{a}^{+}.$ Since $\PPd(\bvsig)$ is concave
on $\mathcal{S}^+_a$, the critical point $\barvsig \in
\mathcal{S}^+_a$ must be a global maximizer of $\PPd(\bvsig)$ on
$\mathcal{S}^+_a$.

On the other hand, by the convexity of $V\left( \bxi\right)$, we
have
\begin{equation}
V\left( \bxi\right) -V\left( \barxi\right) \geq \left\langle
\bxi-\barxi ; \nabla \VV\left( \barxi\right) \right\rangle =\la
 \bxi-\barxi ;  \barvsig  \ra .  \label{Convex}
\end{equation}
Substituting $\boldsymbol{\xi }=\Lambda \left( \bx\right) $ and $
\boldsymbol{\bar{\xi}}=\Lambda \left( \barbx\right)$ into
(\ref{Convex}), we obtain
\begin{equation*}
V\left( \Lambda \left( \bx\right) \right) -V\left( \Lambda \left(
\barbx\right) \right) \geq \la
  \Lambda \left( \bx\right) -\Lambda
\left( \barbx\right) ; \barvsig  \ra.
\end{equation*}
This leads  to
\begin{equation}
\PP \left( \bx\right) -\PP \left( \barbx\right) \geq    \la
\Lambda \left( \bx\right) -\Lambda \left( \barbx\right)  ;
\barvsig \ra + \frac{1}{2} \la \bx , \mathbf{A}\bx \ra -
\frac{1}{2} \la \barbx , \mathbf{A}\barbx \ra -  \la  \bx - \barbx
, \boldsymbol{f}  \ra,  \; \forall \bx \in \real^n . \label{DifPx}
\end{equation}
 By the fact that
\begin{equation}
\boldsymbol{\bar{\varsigma}=}\Lambda \left(
\barbx\right)-\boldsymbol{d},  \label{Kesi}
\end{equation}
  we have
\eb \PP \left( \bx\right) -\PP \left( \barbx\right) \geq
\frac{1}{2} \la \bx, \mathbf{G} \left( \barvsig\right) \bx \ra
-\frac{1}{2} \la \barbx , \mathbf{G} \left( \barvsig\right) \barbx
\ra -  \la  \bx - \barbx , \mathbf{G} \left( \barvsig\right)
\boldsymbol{\bar{x} } \ra .   \label{px} \ee For a fixed $\barvsig
\in \mathcal{S}^+_a$,  the convexity of the complementary gap
function $\boldsymbol{G}_{ap} (\bx, \barvsig) =\frac{1}{2} \la \bx
, \mathbf{G} \left( \barvsig\right) \bx \ra
 $ on $\calX_a$  leads to
\begin{equation}
 \boldsymbol{G}_{ap} (\bx, \barvsig) - \boldsymbol{G}_{ap} (\barbx, \barvsig) \geq \left\langle
\bx-\barbx, \nabla_{\bx} \boldsymbol{G}_{ap}( \barbx , \barvsig )
\right\rangle =  \la \bx-\barbx  , \mathbf{G} \left(
\barvsig\right) \barbx \ra  \;\; \forall \bx \in \real^n.
\label{gx}
\end{equation}
  Therefore,  we have
\begin{equation}
\PP \left( \bx\right) -\PP \left( \barbx\right) \geq \la \bx -
\barbx , \mathbf{G} \left( \barvsig \right) \barbx \ra -  \la \bx
- \barbx, \mathbf{G} \left( \barvsig\right) \boldsymbol{\bar{x} }
\ra =  0     \;\; \forall \bx \in \real^n . \label{Ine3}
\end{equation}
This shows  that $\barbx$ is a global minimizer of Problem $
\left( \mathcal{P} \right)$.  Since it is assumed that
$\barvsig\in \mathcal{S}_a^+$, it follows that (\ref{Global}) is
satisfied.

 We move on to prove the double-min duality statement (\ref{double-min}).

Let $\barvsig$ be a critical point of $\Pi^d(\bvsig)$ and
$\barvsig \in \mathcal{S}_{a}^{-}$. It is easy to verify that
\begin{equation*}
\nabla \PP \left( \bx\right) =\sum\limits_{k=1}^{n}\left(
\frac{1}{2} \bx^{T}\mathbf{B}^{k}\bx-d^{k}\right) \mathbf{B}^{k}
\bx+\mathbf{A}\bx-\boldsymbol{f},
\end{equation*}
\begin{equation}
 \nabla ^{2}\PP \left( \barbx\right) =\mathbf{G} \left(
\barvsig\right) +\mathbf{F} \left( \barbx\right) \mathbf{F} \left(
\barbx\right) ^{T} ,  \label{HesPx}
\end{equation}
where
\begin{equation*}
\mathbf{F} \left( \bx\right) =\left[
\mathbf{B}^{1}\bx,\mathbf{B}^{2}\bx ,\cdots
,\mathbf{B}^{n}\bx\right] .
\end{equation*}

In light of (\ref{GraDual}), $\nabla ^{2}\PP ^{d}\left(
\barvsig\right)$ can be expressed in terms of    $ \barbx=\left[
\mathbf{G} \left( \barvsig\right) \right] ^{-1}\bff$ as follows:
\begin{equation*}
\nabla ^{2}\PP ^{d}\left( \barvsig\right) =-\mathbf{F} \left(
\barbx\right) ^{T}\left[ \mathbf{G} \left( \barvsig \right)
\right] ^{-1}\mathbf{F} \left( \barbx\right) -\mathbf{I},
\end{equation*}
where $\mathbf{I}$ is the identity matrix. If  the critical point
$\barvsig\in \mathcal{S}_{a}^{-}$ is a local minimizer,  we have $
\nabla ^{2}\PP ^{d}\left( \barvsig\right) \succeq 0.$ This leads
to
\begin{equation}
-\mathbf{F} \left( \barbx\right) ^{T}\left[ \mathbf{G} \left(
\barvsig\right) \right] ^{-1}\mathbf{F} \left( \barbx\right)
\succeq \mathbf{I}. \label{Ine4}
\end{equation}
Therefore, $-\mathbf{F} \left( \barbx\right) ^{T}\left[\mathbf{G}
\left( \barvsig\right) \right] ^{-1}\mathbf{F} \left(\barbx\right)
$ is  positive definite
  and  $ \boldsymbol{F} \left( \barbx\right) $ is invertible. By
multiplying $\left( \mathbf{F} \left( \barbx\right) ^{T}\right)
^{-1}$ and $\mathbf{F} \left( \barbx \right)^{-1}$ to the left and
right sides of (\ref{Ine4}), respectively,  we obtain
\begin{equation*}
-\left[ \mathbf{G} \left( \barvsig\right) \right] ^{-1}\succeq
\left( \mathbf{F} \left( \barbx\right) ^{T}\right) ^{-1}  \left(
\mathbf{F} ( \barbx) \right)^{-1}  .
\end{equation*}
According  to Lemma \ref{lem2} in Appendix, the following matrix
inequality is obtained:
\begin{equation*}
\nabla ^{2}\PP \left( \barbx\right) =\mathbf{G} \left(
\barvsig\right) +\mathbf{F} \left( \barbx\right) \mathbf{F} \left(
\barbx\right) ^{T}\succeq 0.
\end{equation*}
By the assumption (\ref{non-deg}), $\barbx =\left[ \mathbf{G}
\left( \barvsig \right) \right] ^{-1}\boldsymbol{f}$ is also a
local minimizer of Problem $ \left( \mathcal{P} \right).$
Therefore,  on a neighborhood $\mathcal{X} _{o}\times
\mathcal{S}_{o}\subset \mathbb{R}^{n}\times \mathcal{S}_{a}^{-}$
of $\left( \barbx,\barvsig\right) ,$ we have
\begin{equation*}
\min_{\bx\in \mathcal{X}_{o}}\PP \left( \bx\right) =\Xi \left(
\barbx,\barvsig\right) =\min_{ \bvsig\in \mathcal{S}_{o}}\PP
^{d}\left( \bvsig\right) .
\end{equation*}
Similarly, we can show that if $ \barbx$ is a local minimizer of
Problem $\left( \mathcal{P}
 \right) ,$  the corresponding $\barvsig$ is also a local
minimizer of Problem $ ( \mathcal{P} ^{d} ) .$

The next task is to prove the double-max duality statement
(\ref{double-max}).

Let $\barvsig\in \mathcal{S}_{a}^{-}$  be a local maximizer of
Problem $ ( \mathcal{P} ^{d} ) $. Then,  we have $ \nabla ^{2}\PP
^{d}\left( \barvsig\right) \preceq 0. $ This gives us
\begin{equation}
\mathbf{F} \left( \barbx\right) ^{T}\left[ \mathbf{G} \left(
\barvsig\right) \right] ^{-1}\mathbf{F} \left( \barbx\right)
+\mathbf{I}\succeq 0 .  \label{MaxPos}
\end{equation}
Now  we have two possible cases regarding the invertibility of
$\mathbf{F} \left( \barbx \right) .$
 If $\mathbf{F} \left(
\barbx\right) $ is invertible, then by  using a similar argument
as presented above, we can show that  the relations
\begin{equation*}
\max_{\bx\in \mathcal{X}_{o}}\PP \left( \bx\right) =\Xi \left(
\barbx,\barvsig\right) =\max_{ \bvsig\in \mathcal{S}_{o}}\PP
^{d}\left( \bvsig\right)
\end{equation*}
hold  on a neighborhood $\mathcal{X}_{o}\times
\mathcal{S}_{o}\subset \mathbb{R} ^{n}\times \mathcal{S}_{a}^{-}$
of $\left( \barbx,\barvsig\right)$.

 If $\mathbf{F} \left( \barbx\right) $ is not
invertible, by   Lemma \ref{lem1}
 in the Appendix, there exists
two orthogonal matrices $\mathbf{E}$ and $\mathbf{K}$ such that
\begin{equation}
\mathbf{F} \left( \barbx\right) =\mathbf{E}\mathbf{D}\mathbf{K},
\label{SVD}
\end{equation}
where
$\mathbf{E}^{T}\mathbf{E}=\mathbf{I}=\mathbf{K}^{T}\mathbf{K} $
and $\mathbf{D} ={\mbox{Diag}}\left( \sigma _{1},\cdots ,\sigma
_{r},0,\cdots ,0\right) $ with  $\sigma _{1}\geq \sigma _{2}\geq
\cdots \geq \sigma _{r}>0 $  and $r={\mbox{rank }}(\mathbf{F}
\left( \barbx\right)) . $ Substituting (\ref{SVD}) into
(\ref{MaxPos}), we obtain
\begin{equation}
-\mathbf{K}^{T}\mathbf{D} ^{T}\mathbf{E}^{T}\left[ \mathbf{G}
\left( \barvsig\right) \right] ^{-1}\mathbf{E}\mathbf{D}
\mathbf{K}-\mathbf{I}\preceq 0. \label{I1}
\end{equation}
Thus,
\begin{equation}
-\mathbf{D} ^{T}\left[ \mathbf{E}^{T}\mathbf{G} \left(
\barvsig\right) \mathbf{E}\right] ^{-1}\mathbf{D}
-\mathbf{I}\preceq 0. \label{I2}
\end{equation}
Applying Lemma \ref{lem4} in Appendix to (\ref{I2}), it follows
that
\begin{equation*}
\mathbf{E}^{T}\mathbf{G} \left( \barvsig\right)
\mathbf{E}+\mathbf{D} \mathbf{D} ^{T}=\mathbf{E}^{T}\mathbf{G}
\left( \barvsig\right)\mathbf{E}+\mathbf{D}
\mathbf{K}\mathbf{K}^{T}\mathbf{D} ^{T}\preceq 0.
\end{equation*}
Finally, we have
\begin{equation*}
\nabla ^{2}\PP \left( \barbx\right) =\mathbf{G} \left(
\barvsig\right) +\mathbf{E}\mathbf{D}
\mathbf{K}\mathbf{K}^{T}\mathbf{D} ^{T}\mathbf{E}^{T}=\mathbf{G}
\left( \barvsig\right) +\mathbf{F} \left(\barbx\right)\mathbf{F}
\left(\barbx\right)^{T}\preceq 0.
\end{equation*}
This means that $\barbx$ is also a local maximizer of Problem $
\left( \mathcal{P} \right) $ under the  assumption
(\ref{non-deg}), i.e., there exists a neighborhood $\mathcal{ \
X}_{o}\times \mathcal{S}_{o}\subset \mathbb{R}^{n}\times
\mathcal{S} _{a}^{-}$ of $\left( \barbx,\boldsymbol{\
\bar{\varsigma}} \right) $ such that
\begin{equation*}
\max_{\bx\in \mathcal{X}_{o}}\PP \left( \bx\right) =\Xi \left(
\barbx,\barvsig\right) =\max_{ \bvsig\in \mathcal{S}_{o}}\PP
^{d}\left( \bvsig\right) .
\end{equation*}

Finally,   we can show, in   a similar way, that if $
\barbx=\left[\mathbf{G} (\barvsig) \right] ^{-1}\boldsymbol{f}$ is
a local maximizer of Problem $\left( \mathcal{P}  \right)$ and
$\barvsig\in \mathcal{S}_{a}^{-}$,
 the corresponding $\barvsig$ is also a local maximizer of
Problem $ ( \mathcal{P} ^{d} ) .$  Therefore, the tri-duality
theorem is proved.  \hfill $\blacksquare$

\begin{remark} The strong triality Theorem \ref{STT} can also be used to identify
saddle points of the primal problem, i.e.
  $\barvsig \in \mathcal{S}_{a}^{-}$  is a saddle point of $\PPd(\bvsig)$ if and only if
    $\barbx = \mathbf{G}(\barvsig)^{-1} \boldsymbol{f} $ is a saddle point of $\PP(\bx)$
   on $\calX_a$.
    Since  the saddle points do not produce computational difficulties in
    numerical optimization, and
     do not exist physically in static systems, these points  are excluded from
    the triality theory.
\end{remark}

\begin{remark}
By the proof of Theorem \ref{STT}, we know that if there exists a
critical point $\barvsig \in \mathcal{S}_{a}^{-}$ such that
$\barvsig$ is a local minimizer of Problem $( \mathcal{P} ^{d} )$,
then $\mathbf{F} (\barbx)$ must be invertible. On the other hand,
if the symmetric  matrices $\{ \boldsymbol{B}^k \}$ are linearly
dependent, then   $\mathbf{F} (\bx)$ is not invertible for any
$\bx \in \mathbb{R}^{n}$. In this case, the corresponding
canonical dual problem $ ( \mathcal{P} ^{d} ) $ has no local
minimizers in $\mathcal{S}_{a}^{-}$, and for any critical point
$\barvsig\in \mathcal{S}_{a}^{-}$, the analytical solution
 $\barbx=\left[ \mathbf{G} \left( \barvsig\right)
\right] ^{-1}\boldsymbol{f}$ is not a local minimizer of
$\PP(\bx)$.
\end{remark}

\section{Refined Triality Theory for General Quartic Polynomial Optimization}

Let us recall the  primal problem and its canonical dual problem
in the general quartic polynomial case ($n\neq m$):
\begin{eqnarray}
\left( \mathcal{P}\right) : &\ext  & \left\{ \PP\left( \bx\right)
=\frac{1}{2}\sum\limits_{k=1}^{m} \left(
\frac{1}{2}\bx^{T}\mathbf{B}^{k}\bx-d^{k}\right)
^{2}+\frac{1}{2}\bx^{T}\mathbf{A}\bx-\bx^{T}\boldsymbol{f} \; | \;
{\bx\in \mathbb{R} ^{n}}\right\}, \;\; \; \;
 \label{Pm} \\
\left( \mathcal{P}^{d}\right) : & \ext  & \left\{ \PP^{d}\left(
\bvsig\right) = -\frac{1}{2}\boldsymbol{f}^{T}\left[
\mathbf{G}\left( \bvsig\right) \right]
^{-1}\boldsymbol{f}-\frac{1}{2}
\bvsig^{T}\bvsig-\bvsig^{T}\boldsymbol{d}\;|\; \bvsig\in
\mathcal{S}_{a}  \subset \real^m  \right\}. \label{Pmd}
\end{eqnarray}
  Suppose that
$\barbx$ and $\barvsig$ are the critical points of Problem $\left(
\mathcal{P}\right) $ and Problem $  ( \mathcal{P}^{d} ) ,$
respectively, where $\barbx = \left[ \mathbf{G}\left(
\barvsig\right) \right] ^{-1} \boldsymbol{f}$. It is easy to
verify that
\begin{equation}
\nabla ^{2}\PP\left( \barbx\right) =\mathbf{G}\left(
\barvsig\right) +\mathbf{F}\left( \barbx\right) \mathbf{F}\left(
\barbx\right) ^{T}    \in \real^{n\times n}  \label{HesPm}
\end{equation}
\begin{equation}
\nabla ^{2}\PP^{d}\left( \barvsig\right) =-\mathbf{F}\left(
\barbx\right) ^{T}\left[ \mathbf{G}\left( \barvsig\right) \right]
^{-1}\mathbf{F}\left( \barbx\right) -\mathbf{I}   \in
\real^{m\times m}. \label{HesPmd}
\end{equation}
In this case,
\begin{equation*}
\mathbf{F}\left( \bx\right) =\left[
\mathbf{B}^{1}\bx,\mathbf{B}^{2} \bx,\cdots
,\mathbf{B}^{m}\bx\right] \in \mathbb{R}^{n \times m} .
\end{equation*}

To continue,   we show the following lemmas.
\begin{lemma} \label{lmwdm}
Suppose that  $m<n$. Let the critical point $\barvsig\in
\mathcal{S}_{a}^{-}$  be a local minimizer of $\PP^{d}(\bvsig)$,
and let $\barbx =  \left[ \mathbf{G}\left( \barvsig\right) \right]
^{-1} \boldsymbol{f}$. Then, there exists a matrix $\mathbf{P}\in
\mathbb{R}^{n \times m}$ with $\rank(\mathbf{P})=m$ such that
 \eb
 \mathbf{P}^{T}\nabla ^{2}\PP\left( \barbx\right)\mathbf{P}
 \succeq \mathbf{0}. \label{WDMP}
  \ee
\end{lemma}
\noindent \textbf{Proof.} By the fact  that  the critical point
$\barvsig\in \mathcal{S}_{a}^{-}$ is a local minimizer of
$\PP^{d}(\bvsig)$, we  have $\nabla \PP^{d}(\barvsig) =0$ and
$\nabla^{2} \PP^{d}(\barvsig) \succeq 0$. It follows that
 \eb
 -\mathbf{F}\left(
\barbx\right) ^{T}\left[ \mathbf{G}\left( \barvsig\right) \right]
^{-1}\mathbf{F}\left( \barbx\right) \succeq \mathbf{I}  \in
\real^{m\times m}. \notag
  \ee
Thus, $\rank(\mathbf{F}(\barbx))=m$. Since $\barvsig\in
\mathcal{S}_{a}^{-}$ and $ \mathbf{F}\left( \barbx\right)
\mathbf{F}\left( \barbx \right) ^{T}\succeq 0,$ there exists a
non-singular matrix $\mathbf{T}\in \mathbb{R} ^{n\times n}$ such
that
\begin{equation}
\mathbf{T}^{T}\mathbf{G}\left( \barvsig\right)
\mathbf{T}={\mbox{Diag }}\left( -\lambda _{1},\cdots ,-\lambda
_{n}\right)
\end{equation}
and
\begin{equation}
\mathbf{T}^{T}\mathbf{F}\left( \barbx\right) \mathbf{F}\left(
\barbx \right) ^{T}\mathbf{T}={\mbox{Diag }}\left( a_{1},\cdots
,a_{m},0,\ldots ,0\right) , \label{SV}
\end{equation}
where $\lambda _{i}>0,$ $i=1,\cdots ,n,$ and $a_{j}>0,$
$j=1,\cdots ,m.$

According to the singular value decomposition theory
\cite{Matrix}, there exist orthogonal matrices $\mathbf{U}$ and
$\mathbf{E}$ such that
\begin{equation*}
\mathbf{T}^{T}\mathbf{F}(\barbx)=\mathbf{U}\left(
\begin{array}{ccc}
  \sqrt{a_1} &  &  \\
   & \ddots &  \\
  &  & \sqrt{a_m} \\
   0 & \cdots & 0  \\
   & \cdots &  \\
    0 & \cdots & 0  \\
\end{array}
\right) \mathbf{E}.
\end{equation*}
Therefore, $\mathbf{U}$ is an identity matrix. Let
\begin{equation*}
\mathbf{R}=\left(
\begin{array}{ccc}
  \sqrt{a_1} &  &  \\
   & \ddots &  \\
  &  & \sqrt{a_m} \\
   0 & \cdots & 0  \\
   & \cdots &  \\
    0 & \cdots & 0  \\
\end{array}
\right) .
\end{equation*}
Then,
\begin{eqnarray*}
\nabla ^{2}\PP^{d}\left( \barvsig\right) &=&-\mathbf{F}\left(
\barbx\right) ^{T}\left[ \mathbf{G}\left( \barvsig\right) \right]
^{-1}\mathbf{F}\left( \barbx\right)
-\mathbf{I}  \\
&=&-\left( \mathbf{F}^{T}\mathbf{T}\right) \left[
\mathbf{T}^{T}\mathbf{G}\left( \barvsig\right) \mathbf{T}\right]
^{-1}\mathbf{T}^{T}\mathbf{F}\left( \barbx
\right) -\mathbf{I} \\
&=&-\mathbf{E}^{T}\mathbf{RU}^{T}\left[\Diag\left( -\lambda
_{1},\cdots ,-\lambda _{n}\right)\right]^{-1}
\mathbf{URE}-\mathbf{I}  \in \real^{m\times m}.
\end{eqnarray*}
Since $\nabla ^{2}\PP^{d}\left( \barvsig\right) \succeq 0$,
$\mathbf{U}$ is an  identity matrix,  and $\mathbf{E}$ is an
orthogonal matrix, we have
\begin{equation*}
-\mathbf{R}[{\mbox{Diag }}\left( \lambda _{1},\cdots ,\lambda
_{n}\right)]^{-1} \mathbf{R}-\mathbf{I}_{m\times m}={\ \mbox{Diag
}}\left( \frac{a_{1}}{\lambda _{1}}-1,\cdots ,\frac{a_{m}}{
\lambda _{m}}-1\right) \succeq 0.
\end{equation*}
Thus, $a_{i}\geq \lambda _{i},$ $i=1,\cdots ,m.$ Note that
 \eb
 \mathbf{T}^{T}\nabla^{2}\PP\left( \barbx\right)
\mathbf{T}={\mbox{Diag }}( a_{1}-\lambda _{1},\cdots
,a_{m}-\lambda_{m},-\lambda _{m+1},\cdots,-\lambda_{n}). \ee Let
$\mathbf{J}=[\mathbf{I}_{m\times m},\mathbf{0}_{m \times
(n-m)}]^{T}$. Then, we have
 \eb
\mathbf{J}^{T}\mathbf{T}^{T}\nabla^{2}\PP\left( \barbx\right)
\mathbf{T}\mathbf{J}={\mbox{Diag }}( a_{1}-\lambda _{1},\cdots
,a_{m}-\lambda_{m}) \succeq \mathbf{0}. \ee
 Let $\mathbf{P}=\mathbf{TJ}$. Clearly, $\rank(\mathbf{P})=m$ and $\mathbf{P}^{T}
 \nabla^{2}\PP\left( \barbx\right) \mathbf{P}= {\mbox{Diag }}( a_{1}-\lambda _{1},\cdots
,a_{m}-\lambda_{m}) \succeq \mathbf{0}$. The proof is completed.
\hfill $\blacksquare$ \medskip

In a similar way, we can prove the following lemma.
\begin{lemma}
Suppose that $m>n$. Let $\barbx=\left[ \mathbf{G}\left(
\barvsig\right) \right] ^{-1}\boldsymbol{f}$ be a critical point,
which is a local minimizer of Problem $\left( \mathcal{P}\right)$,
where $\barvsig \in \mathcal{S}_{a}^{-}$. Then, there exists a
matrix $\mathbf{Q}\in \mathbb{R}^{m \times n}$ with
$\rank(\mathbf{Q})=n$ such that
 \eb
 \mathbf{Q}^{T}\nabla ^{2}\PP^{d}\left( \barvsig\right)\mathbf{Q}
 \succeq \mathbf{0}. \label{WDMQ}
  \ee
\end{lemma}

 Let  $\mathbf{p}_{1},
\cdots, \mathbf{p}_{m}$ be the $m$ column vectors of $\mathbf{P}$
and let $\mathbf{q}_{1}, \cdots, \mathbf{q}_{n}$ be the $n$ column
vectors of $\mathbf{Q}$, respectively. Clearly, $\mathbf{p}_{1},
\cdots, \mathbf{p}_{m}$ are $m$ independent vectors and
$\mathbf{q}_{1}, \cdots, \mathbf{q}_{n}$ are $n$ independent
vectors. We introduce the following two subspaces
 \begin{eqnarray}
  \mathcal{X}_{\flat} &=& \{\bx \in \mathbb{R}^{n}\;|\;\bx = \barbx
  + \theta_{1}\mathbf{p}_{1}+\cdots+\theta_{m}\mathbf{p}_{m},\;
  \theta_{i} \in \mathbb{R},\; i=1,\cdots,m
  \}, \\
  \mathcal{S}_{\flat} & = & \{\bvsig \in \mathbb{R}^{m}\;|\;\bvsig =
  \barvsig + \vartheta_{1}\mathbf{q}_{1}+\cdots+\vartheta_{n}\mathbf{q}_{n},\;
  \vartheta_{i} \in \mathbb{R},\; i=1,\cdots,n
  \} .
  \end{eqnarray}

\begin{theorem}[Refined Triality Theorem] \label{WTT}~~~

Suppose that the assumption (\ref{non-deg}) is satisfied. Let
$\barvsig$ be a critical point of $\PPd(\bvsig)$ and let $\barbx
=\left[ \mathbf{G}\left( \barvsig\right) \right]
^{-1}\boldsymbol{f}$.

If $\barvsig\in \mathcal{S}_{a}^{+},$ then the canonical min-max
duality holds in the strong form:
\begin{equation}
\PP(\barbx) = \min_{\bx\in \mathbb{R} ^{n}}\PP\left(
\bx\right)\Longleftrightarrow \max_{ \boldsymbol{\varsigma }\in
\mathcal{S}_{a}^{+}}\PP^{d}\left( \bvsig \right) =
\PP^{d}(\barvsig) . \label{PmGlobal}
\end{equation}

If $\barvsig\in \mathcal{S}_{a}^{-}$, then there exists a
neighborhood $\mathcal{X}_{o}\times \mathcal{S}_{o}\subset
\mathbb{R} ^{n}\times \mathcal{S}_{a}^{-}$ of $\left(
\barbx,\barvsig\right) $ such that  the double-max duality holds
in the strong form
\begin{equation}
\PP(\barbx) = \max_{\bx\in \mathcal{X}_{o}}\PP\left( \bx\right)
\Longleftrightarrow \max_{ \bvsig\in \mathcal{S}_{o}}\PP^{d}\left(
\bvsig\right) = \PP^{d}\left( \barvsig\right) . \label{Pmdmax}
\end{equation}
However,  the double-min duality statement holds conditionally in
the following  super-symmetrical forms.
\begin{enumerate}
\item   If $m<n$ and $\barvsig\in \mathcal{S}_{a}^{-}$ is a local
minimizer of  $\PPd(\bvsig)$,  then $\barbx =\left[
\mathbf{G}\left( \barvsig\right) \right] ^{-1}\boldsymbol{f}$ is a
saddle point of $\PP(\bx)$  and the double-min duality holds
weakly on $ \mathcal{X}_{o}\cap\mathcal{X}_{\flat} \times
\mathcal{S}_o$, i.e.
 \eb
  \PP(\barbx) = \min_{\bx\in \mathcal{X}_{o}\cap\mathcal{X}_{\flat}}\PP\left( \bx\right)
=\min_{ \bvsig\in \mathcal{S}_{o}}\PP^{d}\left(
\bvsig\right)=\PP^{d}(\barvsig);  \label{WDM1}
 \ee

\item If $m>n$ and $\barbx=\left[ \mathbf{G}\left( \barvsig\right)
\right] ^{-1}\boldsymbol{f}$ is a local minimizer of $\PP(\bx),$
 then $\barvsig$ is a saddle point of $\PPd(\bvsig) $
 and   the double-min duality holds weakly on $\calX_o \times  \mathcal{S}_{o}\cap
\mathcal{S}_{\flat}$, i.e.
 \eb
  \PP(\barbx) = \min_{\bx\in \mathcal{X}_{o}}\PP\left( \bx\right)
=\min_{ \bvsig\in \mathcal{S}_{o}\cap
\mathcal{S}_{\flat}}\PP^{d}\left(
\bvsig\right)=\PP^{d}(\barvsig).\label{WDM2}
 \ee

\end{enumerate}
\end{theorem}
\noindent \textbf{Proof.} The proof of the statements
 (\ref{PmGlobal}) and (\ref{Pmdmax}) are similar to that given for the proof of Theorem
 \ref{STT}. Thus, it suffices to prove the double-min duality
 statements (\ref{WDM1}) and (\ref{WDM2}).

 Firstly,  we
suppose that $m<n$ and $\barvsig$ is a local minimizer of Problem
$(\mathcal{P}^{d})$.
 Define \eb \varphi(t_{1},\cdots,t_{m}) =
\PP(\barbx+t_{1}\barbx_{1}+\cdots+t_{m}\barbx_{m}). \label{phit}
\ee From (\ref {HesPmd}), we obtain
\begin{equation*}
-\mathbf{F}\left( \barbx\right) ^{T}\left[ \mathbf{G}\left(
\barvsig\right) \right] ^{-1}\mathbf{F}\left( \barbx \right)
\succeq \mathbf{I}  \in \real^{m\times m}.
\end{equation*}
Thus, $\mathbf{F}\left( \barbx\right) ^{T}\left[ \mathbf{G}\left(
\barvsig\right) \right] ^{-1}\mathbf{F}\left( \barbx\right) $ is a
non-singular matrix and  ${\mbox{rank }}\left(
\mathbf{F(\barbx)}\right) = m<n$. We claim that $\barbx = \left[
\mathbf{G}\left( \barvsig\right) \right] ^{-1}\boldsymbol{f}$ is
not a local minimizer of Problem $\left( \mathcal{P}\right) .$ On
a contrary, suppose that $\barbx$ is also a local minimizer. Then,
we have
\begin{equation*}
\nabla ^{2}\PP\left( \barbx\right) =\mathbf{G}\left(
\barvsig\right) +\mathbf{F}\left( \barbx\right) \mathbf{F}\left(
\barbx\right) ^{T}\succeq 0.
\end{equation*}
Thus, $ \mathbf{F}\left( \barbx\right) \mathbf{F}\left( \barbx
\right) ^{T}\succeq -\mathbf{G}\left( \barvsig\right) .$ Since
$\barvsig\in \mathcal{S}_{a}^{-}$ and ${ \mbox{rank }}\left(
\mathbf{F}\right) = m$, it is clear that
\begin{equation*}
n={\mbox{rank }}\left( \mathbf{G}\left( \barvsig\right) \right)
={\mbox{rank }}\left( \mathbf{F}\left( \barbx\right)
\mathbf{F}\left( \barbx\right) ^{T}\right) = m.
\end{equation*}
This is a contradiction. Therefore, $\barbx=\left[
\mathbf{G}\left( \barvsig\right) \right] ^{-1}\boldsymbol{f}$ is a
saddle point of Problem $\left( \mathcal{P}\right)$.

It is easy to verify that $\PP(\barbx) = \PP^{d}(\barvsig)$.
 Thus, to prove (\ref{WDM1}), it suffices to prove that  $\mathbf{0}\in
 \mathbb{R}^{m}$ is a local minimizer of the function
 $\varphi(t_{1},\cdots,t_{n})$.

 It is easy to verify that
  \eb  \nabla
\varphi(0,\cdots,0) = [(\nabla
\PP(\barbx))^{T}\mathbf{p}_{1},\cdots,(\nabla
\PP(\barbx))^{T}\mathbf{p}_{m}]^{T}=\nabla
\PP(\barbx))^{T}\mathbf{P}=\mathbf{0}\ee and
 \eb \nabla^{2}\varphi(0,\cdots,0) =
\mathbf{P}^{T}\nabla^{2}\PP\left( \barbx\right) \mathbf{P}. \ee
 In light of Lemma \ref{lmwdm} and the assumption (\ref{non-deg}), it follows that
  $\mathbf{0}\in  \mathbb{R}^{m}$ is, indeed, a local minimizer of the function $\varphi(t_{1},\cdots,
 t_{m})$.

In a similar way, we can establish the case of $m>n$. The proof is
completed. \hfill $\blacksquare$

\begin{remark}
Theorem \ref{WTT} shows that both the canonical min-max and
double-max duality statements hold strongly for general cases; the
double-min duality holds strongly for $n=m$ but weakly for $n\neq
m$ in a symmetrical form. The ``certain additional conditions" are
simply the intersection $\calX_o\bigcap \calX_\flat$ for $m<n$ and
$\calS_o \bigcap \calS_\flat$ for $m>n$. Therefore, the open
problem left in 2003 \cite{gao-opt03,gao-amma03} is solved for the
double-well potential function $\WW(\bx)$.
\end{remark}

 The triality theory has been challenged recently
 in a  series  of more than seven papers, see, for example, \cite{VZ-AA, VZ-JOGO}.
 In the first version of \cite{VZ-JOGO}, Voisei and Zalinescu wrote: ``we consider
that it is important to point out that the main results of this
(triality) theory are false. This is done by providing elementary
counter-examples that lead to think that a correction of this
theory is impossible without falling into trivia". It turns out
that most of these counter-examples simply use the double-well
function $\WW(\bx)$ with $n \neq m $. In fact,
 these counter-examples address the same type of open problem for the double-min duality
 left unaddressed in
\cite{gao-opt03,gao-amma03}. Indeed, by Theorem \ref{WTT}, we know
that both the canonical min-max duality and the double-max duality
hold strongly for the general case $n\neq m$. However, based on
the weak double-min duality, one can easily construct other {\em
V-Z type counterexamples},
 where the strong double-min duality
holds conditionally  when $n\neq m$. Also, interested readers
should find that the references  \cite{gao-opt03,gao-amma03}
never been cited in any one of their papers.

\section{Numerical Experiments}

In this section, some simple numerical examples are presented to
illustrate the canonical duality theory.

\noindent \textbf{Example 1 ($m=n=2$).} Let us first consider
Problem $\left( \mathcal{P} \right) $ with  $n=m=2$.
\begin{equation}
\ext \left\{ \PP\left( \bx\right) =\frac{1}{2}\left[ \left(
\frac{1}{2}\bx^{T}\mathbf{B}^{1}\bx- d_{1}\right) ^{2}+\left(
\frac{1}{2}\bx^{T}\mathbf{B}^{2} \bx-d_{2}\right) ^{2}\right]
+\frac{1}{2}\bx ^{T}\mathbf{A}\bx-\bx^{T}\boldsymbol{f}   | \;
{\bx\in \mathbb{R}^{2}} \right\} , \label{TwoDim}
\end{equation}
where
\begin{equation*}
\mathbf{A}=\left[
\begin{tabular}{ll}
$a_{1}$ & $0$ \\
$0$ & $a_{2}$
\end{tabular}
\right] ,\text{ }\mathbf{B}^{1}=\left[
\begin{tabular}{ll}
$b_{1}$ & $0$ \\
$0$ & $0$
\end{tabular}
\right] ,\text{ }\mathbf{B}^{2}=\left[
\begin{tabular}{ll}
$0$ & $0$ \\
$0$ & $b_{2}$
\end{tabular}
\ \right] ,\text{ }\boldsymbol{f}=\left[ f_{1},f_{2}\right] ^{T}.
\end{equation*}
The canonical dual problem can be expressed as
\begin{equation*}
\PP^{d}\left( \boldsymbol{\varsigma }\right) =-\frac{1}{2}\left(
\frac{ f_{1}^{2}}{a_{1}+\varsigma
_{1}b_{1}}+\frac{f_{2}^{2}}{a_{2}+\varsigma _{2}b_{2}}\right)
-\frac{1}{2}\left( \varsigma _{1}^{2}+\varsigma _{2}^{2}\right)
-\left( d_1 \varsigma _{1}+d_2\varsigma _{2}\right) .
\end{equation*}
Thus,
\begin{equation*}
\nabla \PP^{d}\left( \boldsymbol{\varsigma }\right) =\left[
\begin{tabular}{l}
$\frac{b_{1}f_{1}^{2}}{2\left( a_{1}+\varsigma _{1}b_{1}\right)
^{2}}
-\varsigma _{1}-1$ \\
$\frac{b_{2}f_{2}^{2}}{2\left( a_{2}+\varsigma _{2}b_{2}\right)
^{2}} -\varsigma _{2}-1$
\end{tabular}
\ \right]
\end{equation*}
and
\begin{equation*}
\nabla ^{2}\PP^{d}\left( \boldsymbol{\varsigma }\right) =\left[
\begin{tabular}{ll}
$-b_{1}^{2}f_{1}^{2}\left( a_{1}+\varsigma _{1}b_{1}\right)
^{-3}-1$ & $0$
\\
$0$ & $-b_{2}^{2}f_{2}^{2}\left( a_{2}+\varsigma _{2}b_{2}\right)
^{-3}-1$
\end{tabular}%
\ \ \right] .
\end{equation*}%
Now, we take $b_{1}=b_{2}=f_{1}=f_{2}=1$, $d_{1}=d_{2}=1$ and
$a_{1}=-2,$ $a_{2}=-3$. It is easy to check that $\PP^{d}\left(
\bvsig\right) $ has only one critical point $\barvsig_{1}=\left(
1+\sqrt{2},3.35991198\right) $ in $\mathcal{S}_{a}^{+}$ and four
critical points $ \barvsig_{2}=\left( 3/2,2.59827880\right) ,$  $
\barvsig_{3}=\left( 1-\sqrt{2},-0.45819078\right) ,$  $
\barvsig_{4}=\left( 1-\sqrt{2},2.59827880\right) ,$ $
\barvsig_{5}=\left( 3/2,-0.45819078\right)$ in $S_{a}^{-}$,
respectively. Furthermore, $\barvsig_{2}$ is a local minimizer and
$\barvsig_{3}$ is a local maximizer; the solutions   $
\barvsig_{4}$ and $\barvsig_{5}$ are saddle points of
$\PP^{d}\left( \bvsig\right) $ in $ \mathcal{S}_{a}^{-}. $ Thus,
by Theorem \ref{STT}, we know that $\barbx_{1}=(2.41421356 ,
2.77845711)$ is a global minimizer, while  $\barbx_{2}=(-2 ,
-2.48928859)$ is a local minimizer and $\barbx_{3}=(-0.41421356 ,
$ $   -0.28916855) $ is a local maximizer. The corresponding
values of the cost function are
\begin{equation*}
\PP\left( \barbx_{1}\right) =-14.0421<\PP\left( \barbx_{2}\right)
=-4.3050<\PP\left( \barbx_{3}\right) =0.5971.
\end{equation*}
The graph of $\PP\left( x\right) $ and its contour are depicted in
Figure 1. \medskip

\begin{figure}[tbph]
\centering\includegraphics[width=4.5 in, bb=64 36 537
212]{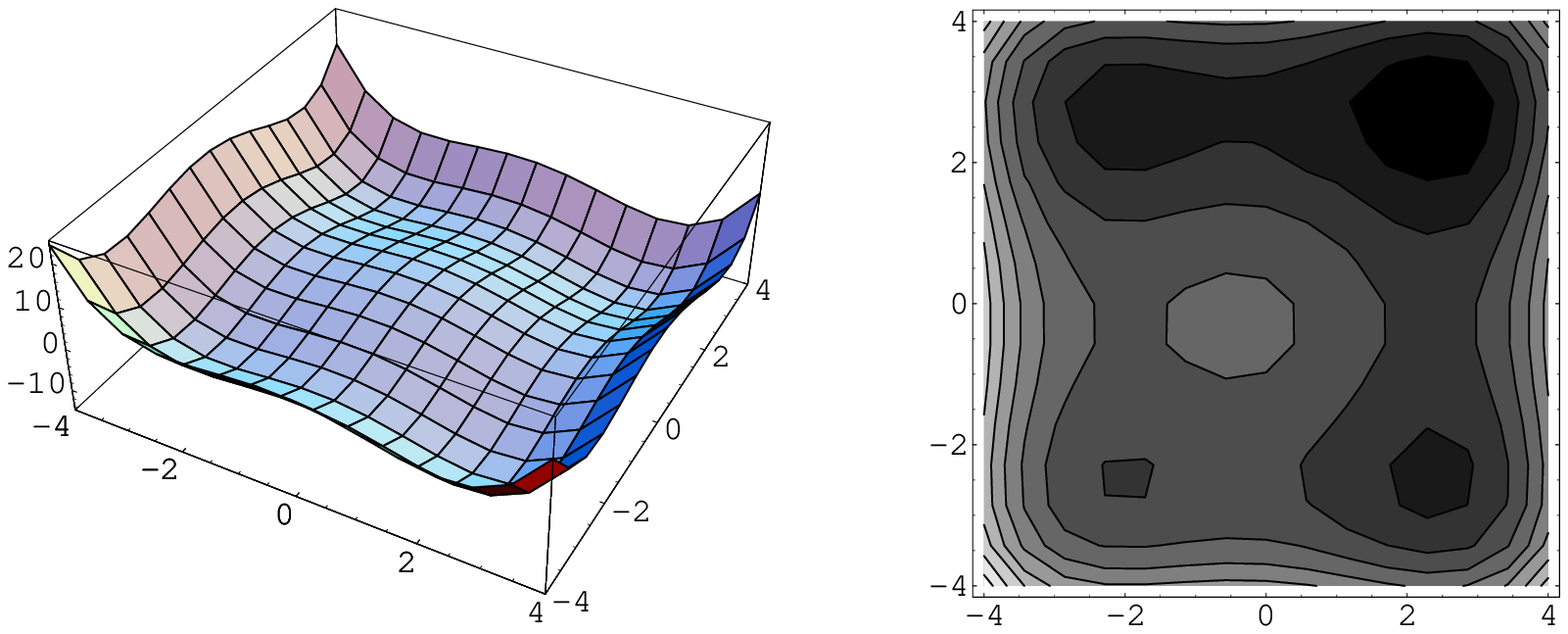}
 \caption{Graph of $\PP(\bx)$ (left) and contours of
$\PP(\bx)$ (right) for Example 1}
\end{figure}

\noindent \textbf{Example 2 ($n=2,m=1$).} We now consider  Problem
$\left( \mathcal{P} \right) $ with $n=2,$ $m=1$, $A = {\mbox{Diag
}}(-0.2,-0.8)$, $B = {\mbox{Diag }}(1,1)$,
$\boldsymbol{f}=(0.9,0.3)^T$, $d=4$. Then, its dual problem is
\begin{equation*}
\PP^{d}\left( \varsigma\right) =-\frac{1}{2}\left( \frac{0.9^{2}
}{-0.2+\varsigma}+\frac{0.3^{2}}{-0.8+\varsigma}\right)
-\frac{1}{2} \varsigma^2 -4 \varsigma.
\end{equation*}
We can verify that $\PP^{d}\left( \varsigma\right) $ has one
critical point $\bar{\varsigma_1} = -0.90489505$ in $S_a^+$ and
two critical points $\bar{\varsigma_2} = -0.12552589$ and
$\bar{\varsigma_3} = -3.974788888$ in $S_a^-$. Furthermore,
$\bar{\varsigma_2}$ is a local minimizer and $\bar{\varsigma_3}$
is a local maximizer of $\PP^{d}(\varsigma)$. According to Theorem
\ref{WTT}, $\barbx_1 = (A + \bar{\varsigma_1} B)^{-1}
\boldsymbol{f} = (1.27678581, 2.86000142)^T$ is the unique global
minimizer of $\PP(\bx)$, $\barbx_2 = (A + \bar{\varsigma_2}
B)^{-1} \boldsymbol{f} = (-2.76475703, -0.32414004)^T$ is a saddle
point and $\barbx_3 = (A + \bar{\varsigma_3} B)^{-1}
\boldsymbol{f}= (-0.21557976, -0.06283003)^T$ is a local maximizer
of $\PP(\bx)$. Let $\boldsymbol{p} = (1,0)^T$ and
$\varphi(\theta)=\PP(\barbx_2+\theta\boldsymbol{p})$. Then, it is
easy to verify that there exists neighborhoods $\mathcal{X}_0
\subset \mathbb{R}$ and $\mathcal{S}_0 \subset \mathbb{R}$ such
that $0 \in \mathcal{X}_0$, $\bar{\varsigma_2}\in \mathcal{S}_0$
and
 \eb
 \min_{\theta\in \mathcal{X}_0} \varphi(\theta) =
 \min_{\varsigma\in \mathcal{S}_0} \PP^d(\varsigma). \notag
 \ee
This example shows that even if $n>m$, the canonical min-max
duality and the double-max duality still hold strongly. However,
the double-min duality statement should be refined into an $m-$
dimensional subspace in this case.

\noindent \textbf{Example 3 ($n=1,m=2$).} We now consider  Problem
$\left( \mathcal{P} \right) $ with $n=1,$ $m=2$, $A = -0.2$, $B^1
= 0.3$, $ B^2 = 0.7$, $d_1 =3$, $d_2=2.7$ and $f=1.4$. Then, its
dual problem is
\begin{equation*}
\PP^{d}\left( \bvsig\right) =-\frac{1}{2}\left( \frac{f^{2}
}{A+\varsigma_1 B^1 + \varsigma_2 B^2}\right) -\frac{1}{2}
(\varsigma_1^2+\varsigma_2^2) -(d_1\varsigma_1+d_2\varsigma_2).
\end{equation*}
We can verify that $\PP^{d}\left( \varsigma\right) $ has one
critical point $\barvsig_1 = (-0.35012607, 3.48303916)^T$ in
$S_a^+$ and two critical points $\barvsig_2 = (-2.98705125,
-2.66978626)^T$ and $\barvsig_3 = (-0.70765026, 2.64881606)^T$ in
$S_a^-$. Furthermore, $\barvsig_2$ is a local maximizer and
$\barvsig_3$ is a saddle point of $\PP^{d}(\bvsig)$. According to
Theorem \ref{WTT}, $\barbx_1 = G(\barvsig_1)^{-1} f = 4.20307342$
is the unique global minimizer of $\PP(\bx)$, $\barbx_2 =
G(\barvsig_2)^{-1} f = -0.29381114$ is a local maximizer of
$\PP(\bx)$. Let $\boldsymbol{q} = (1,0)^T$ and
$\psi(\vartheta)=\PP^d(\barvsig_3+\vartheta \boldsymbol{q})$.
Then, it is easy to verify that there exists neighborhoods
$\mathcal{X}_0 \subset \mathbb{R}$ and $\mathcal{S}_0 \subset
\mathbb{R}$ such that $\barbx_3 \in \mathcal{X}_0$, $0 \in
\mathcal{S}_0$ and
 \eb
 \min_{\bx\in \mathcal{X}_0}\PP(\bx) =
 \min_{\vartheta\in \mathcal{S}_0} \psi(\vartheta). \notag
 \ee
This example shows that if $n<m$, the canonical min-max duality
and the double-max duality still hold strongly. However, the
double-min duality statement should be refined into an $n-$
dimensional subspace in this case.

 \noindent \textbf{Example 4  Linear Perturbation.} Let us consider the following
optimization problem without input ($\boldsymbol{f}= 0 $)
\begin{equation*}
\left( \mathcal{P}_{2}\right) :\ext \left\{ \PP\left( \bx\right)
=\frac{1}{2}\left[ \left( \frac{1 }{2}\left( x_{1}+x_{2}\right)
^{2}-\frac{1}{2}\right) ^{2}+\left( \frac{1}{2} \left(
x_{1}-x_{2}\right) -\frac{1}{2}\right) ^{2}\right] | \; {\bx\in
\mathbb{R} ^{2}} \right\} .
\end{equation*}
Problem $(\mathcal{P}_{2})$ has four global minimizers
$\barbx_{1}=\left( 1,0\right) ,$ $\barbx_{2}=\left( 0,-1\right) ,$
$\barbx_{3}=\left( 0,1\right) ,$ $\barbx_{4}=\left( -1,0\right) $
and the optimal cost value is $0.$ Its canonical dual problem is
\begin{equation*}
\PP^{d}\left( \bvsig\right) =-\frac{1}{2}\left( \varsigma
_{1}^{2}+\varsigma _{2}^{2}\right) -\left( \varsigma
_{1}+\varsigma _{2}\right) .
\end{equation*}
$\PP^{d}\left( \bvsig\right) $ has only one critical point $
\barvsig = \left( -\frac{1}{2},-\frac{1}{2}\right) \in
\mathcal{S}_{a}^{-}$. Furthermore, we can check that $\barbx =
\left[ G\left( \barvsig\right) \right] ^{-1}\boldsymbol{f}= \left(
0,0\right) $ is a local maximizer of Problem $\left( \mathcal{P}
_{2}\right) .$ Thus, we cannot use the canonical dual
transformation method to obtain the global minimizer of Problem
$\left( \mathcal{P}_{2}\right)$ since this problem is in a perfect
symmetrical form without input,
 which allows more than one global minimizer. Now we perturb Problem $\left( \mathcal{P}_{2}\right)
$ as follows.
\begin{equation*}
\left( \mathcal{P}_{2}^{b}\right) :\ext_{\bx\in \real^2}
\PP_{2}\left( \bx\right) =\frac{1}{2}\left[ \left(
\frac{1}{2}\left( x_{1}+x_{2}\right) ^{2}-\frac{1}{2}\right)
^{2}+\left( \frac{1}{2}\left( x_{1}-x_{2}\right)^{2}
-\frac{1}{2}\right) ^{2}\right]  -\left(
x_{1}f_{1}+x_{2}f_{2}\right) .
\end{equation*}%
Its canonical dual function is expressed as
\begin{equation*}
\PP_{2}^{d}\left( \boldsymbol{\varsigma }\right)
=-\frac{1}{8\varsigma _{1}\varsigma _{2}}\left[ \left( \varsigma
_{1}+\varsigma _{2}\right) \left( f_{1}^{2}+f_{2}^{2}\right)
+2\left( \varsigma _{1}-\varsigma _{2}\right) f_{1}f_{2}\right]
-\frac{1}{2}\left( \varsigma _{1}^{2}+\varsigma _{2}^{2}\right)
-\frac{1}{2}\left( \varsigma _{1}+\varsigma _{2}\right).
\end{equation*}
Taking $f_{1}=0.001$, $f_{2}=0.005$ and solving $\nabla
\PP_{2}^{d}\left( \barvsig\right) =0, $ the results obtained are
listed in Table 1. We can see that $\barvsig=\left(
0.00299107,0.00199602\right) \in \mathcal{S}_{a}^{+}$ and
$\PP_{2}^{d}\left( \bvsig\right) =-0.00500648$. Thus, $\barbx =
\left[ \mathbf{G}\left( \barvsig\right) \right] ^{-1}
\boldsymbol{f} = \left( 0.000495793,1.00249\right) $ is the global
minimizer of Problem $ \left( \mathcal{P}_{2}^{b}\right) $.
Clearly, this $\barbx$ is very close to $\barbx_{3}$. If we take
$f_{1}=0.001$, $ f_{2}=-0.005,$ the global minimizer of Problem
$\left( \mathcal{P} _{2}^{b}\right) $ is $\barbx = \left[
\mathbf{G}\left( \barvsig\right) \right] ^{-1}\boldsymbol{f}=
\left( 0.000496288,-1.00249\right) $ which is close to $\barbx
_{2}=\left( 0,-1\right) $. This example shows that if the
canonical dual problem has no critical point in
$\mathcal{S}_{a}^{+},$ a linear  perturbation could be used to
solve the primal problem.
\begin{table}[tbph]
\caption{Numerical results for Example 3}\centering
\includegraphics[width=5 in,bb=52 520 551 744]
{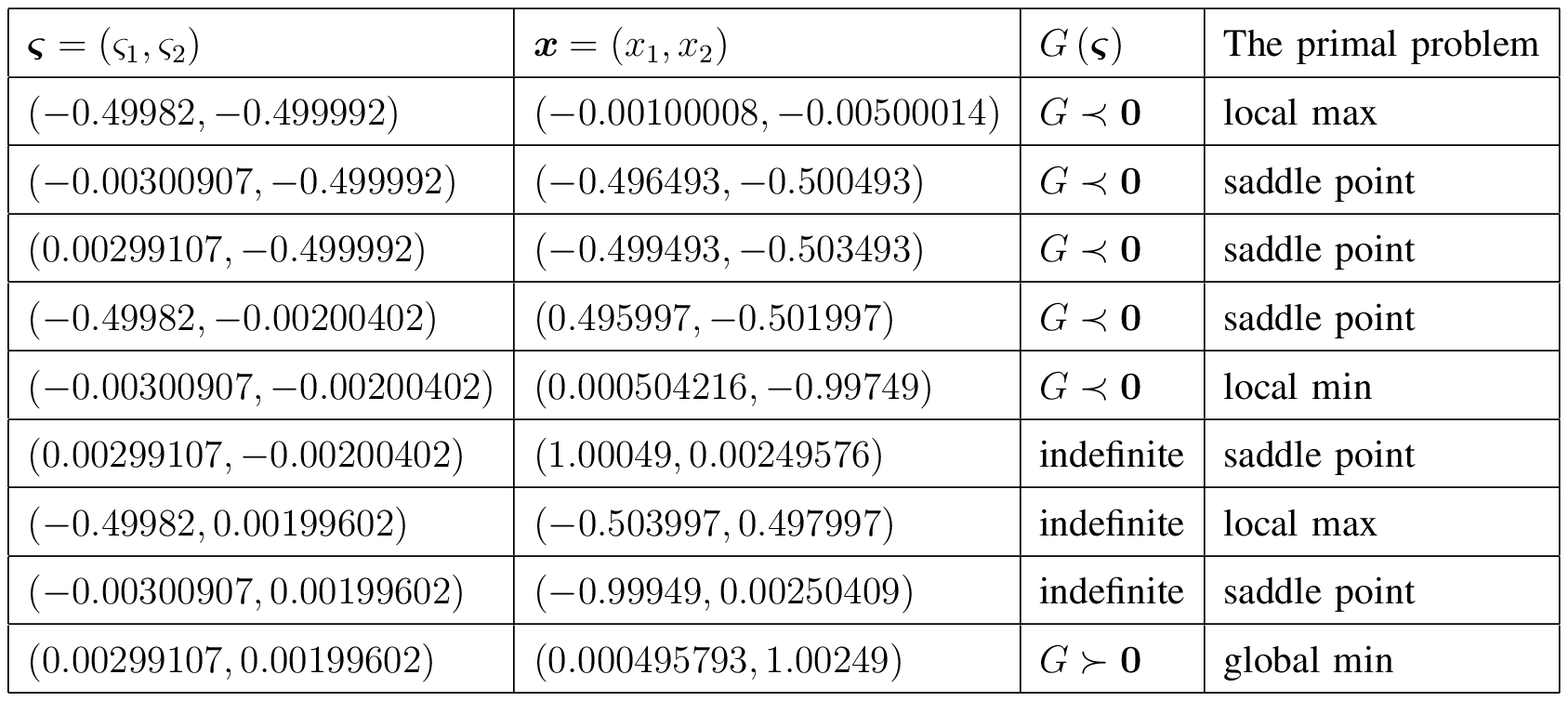}
\end{table}

\section{Appendix}

In this Appendix, we present several lemmas which are needed for
the proofs of Theorem \ref{STT} and Theorem \ref{WTT}.

\begin{lemma} [Singular value decomposition  \cite{Matrix}] \label{lem1}
For any given $\mathbf{G} \in \mathbb{R} ^{n\times n}$ with $
\rank(\mathbf{G}) = r, $
 there exist $ \mathbf{U} \in
\mathbb{R} ^{n\times n},$ $\mathbf{D} \in \mathbb{R} ^{n\times n}$
and $\mathbf{R} \in \mathbb{R} ^{n\times n}$ such that
\begin{equation*}
\mathbf{G} = \mathbf{U}\mathbf{D} \mathbf{R} ,
\end{equation*}
where $\mathbf{U}$,  $\mathbf{R}$ are orthogonal matrices, i.e.,
$\mathbf{U}^{T}\mathbf{U}=\mathbf{I}=\mathbf{R}^{T}\mathbf{R}$,
and $ \mathbf{D} = {\mbox{Diag }}( \sigma _{1},\cdots ,\sigma
_{r},0,\cdots ,0) ,$ $\sigma _{1}\geq \sigma _{2}\geq \cdots \geq
\sigma _{r}>0$ .
\end{lemma}

\begin{lemma} \label{lem2} Suppose that $\mathbf{G} $ and $\mathbf{U}$ are positive
definite. Then, $\mathbf{G} \succeq \mathbf{U} $ if and only if
$\mathbf{U}^{-1}\succeq \mathbf{G}^{-1}$.
\end{lemma}

\noindent \textbf{Proof.} The proof is trivial and is  omitted
here.

\begin{lemma} [Proposition 2.1 in \cite{Schur}] \label{lem3}
For any given symmetric matrix $\mathbf{M}$  expressed in the form
\begin{equation*}
\mathbf{M}=\left[
\begin{tabular}{ll}
$\mathbf{M}_{11}$ & $\mathbf{M}_{12}$ \\
$\mathbf{M}_{21}$ & $\mathbf{M}_{22}$%
\end{tabular}
\right]
\end{equation*}
such that  $\mathbf{M}_{22}\succ 0.$ Then, $\mathbf{M}\succeq 0$
if and only if $
\mathbf{M}_{11}-\mathbf{M}_{12}\mathbf{M}_{22}^{-1}\mathbf{M}_{21}\succeq
0.$
\end{lemma}

The following lemma plays a key role in the proof of Theorem
\ref{STT} and Theorem \ref{WTT}.

\begin{lemma} \label{lem4}Suppose that $\mathbf{P}\in \real^{n\times n},
\; \mathbf{U}\in \real^{m\times m} , $ and $\mathbf{D}\in
\real^{n\times m}$. Furthermore,
\[
\mathbf{D} =\left[
\begin{tabular}{ll}
$\mathbf{D} _{11}$ & $\mathbf{0}_{r\times \left( m-r\right) }$ \\
$\mathbf{0}_{\left( n-r\right) \times r}$ & $\mathbf{0}_{\left(
n-r\right) \times \left(
m-r\right) }$%
\end{tabular}%
\ \right] \in \mathbb{R}^{n\times n} ,
\]
  where  $\mathbf{D}
_{11}\in \mathbb{R} ^{r\times r}$ is nonsingular,
$r=\rank(\mathbf{D})$, and
 \[
\mathbf{P}=\left[
\begin{tabular}{ll}
$\mathbf{P}_{11}$ & $\mathbf{P}_{12}$ \\
$\mathbf{P}_{21}$ & $\mathbf{P}_{22}$
\end{tabular}%
\right]  \prec 0, \;\;\; \mathbf{U}=\left[
\begin{tabular}{ll}
$\mathbf{U}_{11}$ & $\mathbf{0}_{r\times \left( m-r\right) }$ \\
$\mathbf{0}_{(m-r)\times r }$ & $\mathbf{U}_{22}$
\end{tabular}
\ \right] \succ 0,
\]
$\mathbf{P}_{ij}$ and $\mathbf{U}_{ii}$, $i,j=1,2$, are of
appropriate dimension matrices.  Then,
\begin{equation}
\mathbf{P}+\mathbf{D} \mathbf{U}\mathbf{D} ^{T}\preceq
0\Longleftrightarrow -\mathbf{D} ^{T}\mathbf{P}^{-1}\mathbf{D}
-\mathbf{U}^{-1}\preceq 0. \label{Lem4Ine}
\end{equation}
\end{lemma}

\noindent \textbf{Proof.} Suppose that $\mathbf{P}+\mathbf{D}
\mathbf{UD} ^{T}\preceq 0.$ Then,
\begin{equation*}
-\mathbf{P}-\mathbf{D UD} ^{T}=\left[
\begin{tabular}{ll}
$-\mathbf{P}_{11}-\mathbf{D} _{11}\mathbf{U}_{11}\mathbf{D}
_{11}^{T}$ & $-\mathbf{P}_{12}$ \\
$-\mathbf{P}_{21}$ & $-\mathbf{P}_{22}$%
\end{tabular}
\right] \succeq 0.
\end{equation*}
Since $\mathbf{P}=\left[
\begin{tabular}{ll}
$\mathbf{P}_{11}$ & $\mathbf{P}_{12}$ \\
$\mathbf{P}_{21}$ & $\mathbf{P}_{22}$
\end{tabular}
\right] \prec 0,$ it follows that $-\mathbf{P}_{22}\succ 0.$ By
Lemma \ref{lem3}, we have the following inequality
\begin{equation}
-\mathbf{P}_{11}-\mathbf{D} _{11}\mathbf{U}_{11}\mathbf{D}
_{11}^{T}+\mathbf{P}_{12}\mathbf{P}_{22}^{-1}\mathbf{P}_{21}\succeq
0  \label{Lem3ine1}
\end{equation}
which leads to
\begin{equation*}
-\mathbf{P}_{11}+\mathbf{P}_{12}\mathbf{P}_{22}^{-1}\mathbf{P}_{21}\succeq
\mathbf{D} _{11}\mathbf{U}_{11}\mathbf{D} _{11}^{T}.
\end{equation*}
Since $\mathbf{P}\prec 0$ and $\mathbf{U}\succ 0,$ it follows from
Lemma \ref{lem2} that
\begin{equation*}
\left(
-\mathbf{P}_{11}+\mathbf{P}_{12}\mathbf{P}_{22}^{-1}\mathbf{P}_{21}\right)
^{-1}\preceq \left(\mathbf{D} _{11}^{T}\right)
^{-1}\mathbf{U}_{11}^{-1}\mathbf{D} _{11}^{-1}.
\end{equation*}
Thus,
\begin{equation}
\mathbf{D} _{11}^{T}\left(
-\mathbf{P}_{11}+\mathbf{P}_{12}\mathbf{P}_{22}^{-1}\mathbf{P}_{21}\right)
^{-1}\mathbf{D} _{11}\preceq \mathbf{U}_{11}^{-1}.
\label{Lem3ine2}
\end{equation}
Note that
\begin{equation*}
\mathbf{P}^{-1}=\left[
\begin{tabular}{ll}
$\left(
\mathbf{P}_{11}-\mathbf{P}_{12}\mathbf{P}_{22}^{-1}\mathbf{P}_{21}\right)
^{-1}$ & $\mathbf{P}_{11}^{-1}\mathbf{P}_{12}
\left( \mathbf{P}_{21}\mathbf{P}_{11}^{-1}\mathbf{P}_{12}-\mathbf{P}_{22}\right) ^{-1}$ \\
$\left(
\mathbf{P}_{21}\mathbf{P}_{11}^{-1}\mathbf{P}_{12}-\newline
\mathbf{P}_{22}\right) ^{-1}\mathbf{P}_{21}\mathbf{P}_{11}^{-1}$ &
$\left(
\mathbf{P}_{22}-\mathbf{P}_{21}\mathbf{P}_{11}^{-1}\mathbf{P}_{12}\right)
$
\end{tabular}
\right] .
\end{equation*}
By virtue of Lemma \ref{lem3}, we obtain
\begin{equation*}
-\left[
\begin{tabular}{ll}
$\mathbf{D} _{11}$ & $\mathbf{0}$ \\
$\mathbf{0}$ & $\mathbf{0}$
\end{tabular}
\right] ^{T}\mathbf{P}^{-1}\left[
\begin{tabular}{ll}
$\mathbf{D} _{11}$ & $\mathbf{0}$ \\
$\mathbf{0}$ & $\mathbf{0}$%
\end{tabular}
\right] \preceq \left[
\begin{tabular}{ll}
$\mathbf{U}_{11}^{-1}$ & $\mathbf{0}$ \\
$\mathbf{0}$ & $\mathbf{U}_{22}^{-1}$%
\end{tabular}
\right] =\mathbf{U}^{-1},
\end{equation*}
i.e., the right hand side of (\ref{Lem4Ine}) holds.

 In a similar
way, we can show that if $-\mathbf{D}
^{T}\mathbf{P}^{-1}\mathbf{D} -\mathbf{U}^{-1}\preceq 0,$ then
$\mathbf{P}+\mathbf{DU}\mathbf{D} ^{T}\preceq 0.$ The
proof is thus completed.  \hfill $\blacksquare$ \\

\section{Conclusion Remarks}
In this paper, we presented a rigorous proof of the double-min
duality in the triality theory for a quartic polynomial
optimization problem based on elementary linear algebra. Our
results show that under some proper assumptions, the triality
theory for a class of quartic polynomial optimization problems
holds strongly in the tri-duality form if the primal problem and
its canonical dual have the same dimension. Otherwise, both the
canonical min-max and the double-max   still hold strongly, but
the double-min duality holds weakly in a symmetric form.

\section*{Acknowledgments} The main results of this paper were announced in the 2nd World
Congress  of Global Optimization, July 3-7, 2011, Chania, Greece.
The authors are sincerely indebted to  Professor Hanif Sherali at
Virginia Tech for his valuable comments and suggestions.
 David Gao's research is supported by US Air Force
Office of Scientific Research under the grant AFOSR
FA9550-10-1-0487.
 Changzhi Wu was supported by National Natural
Science Foundation of China under the grant \# 11001288, the Key
Project of Chinese Ministry of Education under the  grant \#
210179, SRF for ROCS, SEM, Natural Science Foundation Project of
CQ CSTC under the  grant \# 2009BB3057 and CMEC under the grant \#
KJ090802.

\end{document}